\documentclass[11pt]{article}
\usepackage[english]{babel}
\usepackage{times}
\usepackage[T1]{fontenc}
\usepackage{amsmath, amsthm, amsfonts}
\usepackage{amssymb}
\usepackage{amscd}
\usepackage{verbatim}

\title{Sharp Trace  Hardy-Sobolev-Maz'ya Inequalities \\ and the Fractional Laplacian}

\author{\Large Stathis Filippas$^{1,4}$,~~  Luisa Moschini$^{2}$~
\&~ Achilles Tertikas$^{3,4}$  \\
                                                                           \\
        Department of Applied Mathematics$^{1}$ \\
         University of Crete,
         71409 Heraklion,  Greece \\
        filippas@tem.uoc.gr\\
                                          \\
                                          Dipartimento di Scienze di Base  ed Applicate per l'Ingegneria $^2$,\\
 University of Rome ``La Sapienza'', 00185 Rome, Italy \\
moschini@dmmm.uniroma1.it \\
\\
 Department of Mathematics$^{3}$ \\
         University of Crete,
         71409 Heraklion,  Greece \\
          tertikas@math.uoc.gr\\
                                     \\
        Institute of Applied and Computational Mathematics$^4$, \\
        FORTH, 71110 Heraklion, Greece \\
    \\ }
\begin{document}
\date{}
\newcommand{\ana}{\nabla}
\newcommand{\R}{I \!  \! R}
\newcommand{\N}{I \!  \! N}
\newcommand{\Ren}{ I \! \! R^n}
\maketitle

%Equations etc

\newcommand{\be}{\begin{equation}}
\newcommand{\ee}{\end{equation}}
\newcommand{\bea}{\begin{eqnarray}}
\newcommand{\eea}{\end{eqnarray}}
\newcommand{\bean}{\begin{eqnarray*}}
\newcommand{\eean}{\end{eqnarray*}}
\newcommand{\la}{\label}
\newcommand{\bx}{\bar{x}}

%Greek Letters

\newcommand{\xa}{\alpha}
\newcommand{\xb}{\beta}
\newcommand{\xg}{\gamma}
\newcommand{\xG}{\Gamma}
\newcommand{\xd}{\delta}
\newcommand{\xD}{\Delta}
\newcommand{\xe}{\varepsilon}
\newcommand{\xz}{\zeta}
\newcommand{\xh}{\eta}
\newcommand{\Th}{\Theta}
\newcommand{\xk}{\kappa}
\newcommand{\xl}{\lambda}
\newcommand{\xL}{\Lambda}
\newcommand{\CC}{\mathcal{C}}
\newcommand{\xm}{\mu}
\newcommand{\xn}{\nu}
\newcommand{\ks}{\xi}
\newcommand{\KS}{\Xi}
\newcommand{\xp}{\pi}
\newcommand{\xP}{\Pi}
\newcommand{\xr}{\rho}
\newcommand{\xs}{\sigma}
\newcommand{\xS}{\Sigma}
\newcommand{\xf}{\phi}
\newcommand{\xF}{\Phi}
\newcommand{\ps}{\psi}
\newcommand{\PS}{\Psi}
\newcommand{\xo}{\omega}
\newcommand{\xO}{\Omega}
\newcommand{\tA}{\tilde{A}}
\newcommand{\tC}{\bar{k}_s}
\newcommand{\finedim}{{\hfill $\Box$}}

%integrals
\newcommand{\dino}{\int_0^{+\infty}\int_{\xO}}
\newcommand{\dinoC}{\int_0^{+\infty}\int_{\CC \xO}}
\newcommand{\dinoc}{\int_0^{+\infty}\int_{\xO^c}}
\newcommand{\dinr}{\int_0^{+\infty}\int_{\Ren}}
\newcommand{\dd}{dxdy}
\newcommand{\ra}{\rightarrow}
\newcommand{\rft}{\rightarrow +\infty}
\newcounter{newsection}

%Theorems Lemmas etc

\newtheorem{theorem}{Theorem}[section]
\newtheorem{lemma}[theorem]{Lemma}
\newtheorem{prop}[theorem]{Proposition}
\newtheorem{coro}[theorem]{Corollary}
\newtheorem{defin}[theorem]{Definition}
\newcounter{newsec} \renewcommand{\theequation}{\thesection.\arabic{equation}}

\begin{abstract}
In this work we  establish  trace Hardy and trace
Hardy-Sobolev-Maz'ya inequalities with best Hardy constants,
 for domains satisfying
suitable geometric assumptions such as mean convexity or
convexity. We then use them to produce fractional
Hardy-Sobolev-Maz'ya inequalities with best Hardy  constants for
various fractional Laplacians. In the case where the domain is the
half space our results cover the full range of the exponent $s \in
(0,1)$ of the fractional  Laplacians. We answer in particular an
open problem raised by Frank and Seiringer \cite{FS}.
\end{abstract}



\tableofcontents

\setcounter{equation}{0}
\section{Introduction and Main Results}\la{sec1}

The Hardy inequality in the upper half space asserts that
\begin{equation}
\int_{\Ren_+}|\nabla u|^2dx \geq \frac{1}{4}
\int_{\Ren_+}\frac{|u|^2}{x_n^2}dx , \quad \quad \quad u
\in C^{\infty}_0(\Ren_+), \label{maz}
\end{equation}
where $\Ren_+ = \{ (x_1, \ldots, x_n): x_n >0 \}$
 denotes the upper half-space, and $\frac14$ is the best possible constant.

If $\xO \subset \Ren$    and $d(x) = {\rm dist}(x, \partial \xO)$
then there are two main directions towards establishing Hardy
inequalities. One direction is to find proper regularity
assumptions on the boundary of $\xO$ that imply
 the existence of a positive constant
$C_{\xO}$ such that
\[
\int_{\xO}|\nabla u|^2dx \geq C_{\xO}
\int_{\xO}\frac{|u|^2}{d^2(x)}dx \ , ~~~~~~~~~~~~~~~~~~  u \in C^{\infty}_{0}(\xO)  \ .
\]
In this direction we refer to \cite{A}, \cite{KK} and references therein.

A second  direction aims at finding geometric assumptions on $\xO$ that imply the Hardy inequality
with best constant  $ \frac14$, that is
\be\la{hardy}
\int_{\xO}|\nabla u|^2dx \geq  \frac14
\int_{\xO}\frac{|u|^2}{d^2(x)}dx \ , ~~~~~~~~~~~~~~~~~~  u \in C^{\infty}_{0}(\xO)  \ .
\ee
The standard geometric assumption here is convexity  of $\xO$, see, e.g., \cite{D1}, \cite{D2}, \cite{BM}.
However inequality  (\ref{hardy})  remains true
 under the weaker assumption
\be\la{c}
- \Delta d(x) \geq 0, ~~~~~~~~~~~~~ x \in \xO  \ .
\ee
This is meant in the distributional sense. We refer to  \cite{BFT} where this condition arises in a natural way.
 In fact  condition (\ref{c}) is equivalent to convexity in two space dimensions,  but it is
weaker than convexity  for $n \geq 3$, since  any convex domain satisfies  (\ref{c}) whereas there are
nonconvex domains that satisfy (\ref{c})  \cite{AK}. We emphasize that there is  no need for  further  regularity
assumptions on $\xO$.  In case $\partial \xO$ is   $C^2$,  condition  (\ref{c}) is recently shown to be equivalent
to the mean convexity of $\partial \xO$, that is
 $(n-1)H(x) = - \Delta d(x) \geq 0$ for $x \in \partial \xO$, see \cite{LLL},
\cite{P}.

If in addition to (\ref{c}) the domain  $\xO$  is a  $C^2$ domain with finite inner radius
 then it has been established
 that  one can combine
the Sobolev and the  Hardy inequality, the latter with best constant. More precisely, for $n \geq 3$
 there exists a positive
constant $c$ such that
\begin{equation}
\int_{\xO}|\nabla u|^2dx \geq \frac{1}{4}
\int_{\xO}\frac{|u|^2}{d^2(x)}dx + c \left(
\int_{\xO}|u|^{\frac{2n}{n-2}}dx \right)^{\frac{n-2}{n}}, \quad u
\in C^{\infty}_0(\xO)  \ , \label{maz1}
\end{equation}
see \cite{FMT2}. In \cite{Gk} Hardy-Sobolev-Maz'ya inequalities
are established under a different geometric  assumption than
(\ref{c}), that allows infinite inner radius.  Frank and Loss
established in   \cite{FL}  inequality (\ref{maz1}) with a
constant $c$ independent of $\xO$, when $\xO$ is convex.

Recently, a lot of attention is attracted by   the fractional Laplacian. For $s \in (0,1)$ it is defined as
follows
\be\la{1.2}
(-\Delta)^s f(x)=c_{n,s}~ P.V.~ \int_{\R^n}
\frac{f(x)-f(\xi)}{|x-\xi|^{n+2s}} d\xi \ ,
\ee
 where P.V. stands for the Cauchy
principal value and
\be\la{1.cns}
c_{n,s}=\frac{s 2^{2s}\Gamma
\left(\frac{n+2s}{2}\right)}{\Gamma(1-s) \pi^{\frac{n}{2}}}  \ .
\ee

There are other ways for defining the fractional Laplacian, as for instance via the Fourier transform.
We note that the fractional Laplacian is a non local operator and this raises several technical difficulties.
However, there is a way of studying various properties of the  fractional Laplacian via the Dirichlet to
Neumann map. This  has been recently studied by Caffarelli and Silvestre \cite{CS}, and it will be central in this work.
Let us briefly recall the approach in \cite{CS}, where by
adding a new variable $y$,  they relate  the fractional Laplacian to a local operator. For
any function $f$  one solves the following  extension  problem
\bea\la{ext1}
 div(y^{1-2s} \nabla_{(x,y)} u(x,y)) &  =  &  0,~~  \ \
\R^n\times (0,\infty), \\ u(x,0) &  =  &  f(x),~~ \ \ \R^n  \la{ext2} \ ,
\eea
the natural energy of which is given by
\[
J[u]=  \dinr y^{1-2s} |\ana_{(x,y)} u(x,y)|^2 \dd.
\]
Then, up to a normalizing factor $C$ one establishes that
\[
-\lim_{y\to 0^{+}} y^{1-2s} u_y(x,y) = C(-\Delta)^s f(x) \ .
\]

Our interest in this work  is to  study  the fractional Laplacian
defined in subsets of $\Ren$ and in particular to   establish
Hardy and  Hardy-Sobolev-Maz'ya inequalities there. There is a lot
of interest in  fractional Laplacian in  subsets of $\Ren$ coming
from various applications, as for instance  censored stable
processes  and killed stable processes \cite{CSo}, \cite{BBC},
\cite{CKS1}, \cite{CKS2}, Gamma convergence and phase transition problems
\cite{ABS}, \cite{G}, \cite{SV1}, \cite{SV2}, \cite{PSV} and
nonlinear PDE theory \cite{CT}, \cite{T}, \cite{CC}. In \cite{BD}
it was conjectured that the best Hardy constant in the  case of
the fractional Laplacian associated to  a  censored stable process
is the same for all convex domains. In \cite{FS} it was posed  the
question establishing fractional Hardy-Sobolev-Maz'ya inequalities
for the half space.

Contrary to the case of the full space $\Ren$, there are several different
  fractional Laplacians  that one can define
 on a domain $\xO \subsetneqq  \Ren$. In particular in the above mentioned references
three different fractional Laplacians appear.
 In all cases we will use
 the Dirichlet to  Neumann map after identifying  the proper extension problem. Throughout this work
we assume that the domain $\xO$ is a uniformly Lipschitz domain; for the precise definition see Section 2.

We start with  the fractional Laplacian that appears in \cite{CT}, \cite{T}, \cite{CC}. The proper
extension problem  in this case is to  consider test functions in  $C^{\infty}_{0}(\xO \times \R) $.
At this point we recall that the inner radius of a domain $\xO$ is defined as $R_{in} := \sup _{x \in \xO} d(x)$.
We say that the  domain  $\xO$ has finite inner radius whenever  $R_{in} < \infty$.
Our  first result concerns the extended problem and reads:

\begin{theorem}\la{th11}({\bf Trace Hardy \& Trace Hardy-Sobolev-Maz'ya I}) \\
\noindent  Let $\frac12 \leq s <1$, $n \geq 2$  and $\xO \subsetneqq \Ren$  be a domain.  \\
\noindent  {\bf (i)} If in  addition  $\xO$ is
 such that
\be\la{c2}
- \Delta d(x) \geq 0, ~~~~~~~~~~~~~ x \in \xO  \ ,
\ee
  then for all   $u \in C^{\infty}_{0}(\xO \times \R) $ there holds
\be\la{1.A1}
 \dino y^{1-2s} |\ana_{(x,y)} u(x,y)|^2 \dd  \geq \bar{d}_{s}  \int_{\xO} \frac{u^2(x,0)}{d^{2s}(x)}dx \ ,
\ee
with
\be\la{1.A2}
 \bar{d}_{s} :=  \frac{2 \Gamma \left(1-s \right) \Gamma^2 \left(\frac{3+2s}{4} \right)}
{\Gamma^2 \left(\frac{3-2s}{4} \right) \Gamma \left( s \right)} \ .
\ee
\noindent  {\bf(ii)}
 Suppose there exists a point $x_0 \in \partial \xO$ and $r>0$ such that the
part of the boundary $\partial \xO  \cap B(x_0,r)$ is $C^1$ regular. Then
\[
 \bar{d}_{s} \geq
 \inf_{u \in C^{\infty}_{0}(\xO \times \R)}
 \frac{\dino y^{1-2s} |\ana u|^2 \dd}{\int_{\xO} \frac{u^2(x,0)}{d^{2s}(x)} dx}  \ .
\]
In particular $ \bar{d}_{s}$ in (\ref{1.A1}) is the best constant. \\
\noindent  {\bf(iii)} If $\xO$ is a uniformly Lipschitz  domain with finite inner radius  satisfying (\ref{c2}),
and $s \in (\frac12,1)$,
then there exists a positive constant $c$  such that
for all   $u \in C^{\infty}_{0}(\xO \times \R) $ there holds
\be\la{1.A3}
 \dino y^{1-2s} |\ana_{(x,y)} u(x,y)|^2 \dd  \geq  \bar{d}_{s}
 \int_{\xO} \frac{u^2(x,0)}{d^{2s}(x)}dx + c
 \left(  \int_{\xO} | u(x,0)|^{\frac{2n}{n-2s}} dx \right)^{\frac{n-2s}{n}}  \ .
\ee
\end{theorem}
% \vspace{2mm}

Actually, in the case of half space $\xO= \R^n_{+}$ we  establish  a much stronger result covering the full range
$s \in (0,1)$. In particular we have

%\vspace{2mm}

\begin{theorem}\la{th12}({\bf Half Space,  Trace  Hardy-Sobolev-Maz'ya I}) \\
Let $0<s<1$ and  $n \geq 2$. \\
\noindent {\bf (i)} For all   $u \in C^{\infty}_{0}(\R^n_{+} \times \R) $ there holds\
\be\la{1.13}
 \int_{0}^{\infty} \int_{\R^n_{+}}  y^{1-2s} |\ana_{(x,y)} u(x,y)|^2 \dd  \geq  \bar{d}_{s}
 \int_{\R^n_{+}} \frac{u^2(x,0)}{x_n^{2s}}dx \ ,
\ee
with
\be\la{1.13d}
 \bar{d}_{s} :=  \frac{2 \Gamma \left(1-s \right) \Gamma^2 \left(\frac{3+2s}{4} \right)}
{\Gamma^2 \left(\frac{3-2s}{4} \right) \Gamma \left( s \right)} \ .
\ee
\noindent {\bf (ii)} The constant $\bar{d}_{s}$ in (\ref{1.13}) is sharp, that is
\[
 \bar{d}_{s} =
 \inf_{u \in C^{\infty}_{0}(\R^n_{+} \times \R)}
 \frac{ \int_{0}^{\infty} \int_{\R^n_{+}}  y^{1-2s} |\ana u|^2 \dd}{\int_{\R^n_{+}} \frac{u^2(x,0)}{x_n^{2s}} dx}  \ .
\]
\noindent {\bf (iii)} There exists a positive constant $c$ such that
 for all   $u \in C^{\infty}_{0}(\R^n_{+} \times \R) $ there holds
\be\la{1.13sob}
 \int_{0}^{\infty} \int_{\R^n_{+}}  y^{1-2s} |\ana_{(x,y)} u(x,y)|^2 \dd  \geq  \bar{d}_{s}
 \int_{\R^n_{+}} \frac{u^2(x,0)}{x_n^{2s}}dx + c
 \left(  \int_{\R^n_{+}} | u(x,0)|^{\frac{2n}{n-2s}} dx \right)^{\frac{n-2s}{n}}  \ .
\ee
\end{theorem}

We  will apply Theorem \ref{th11}  to
 the   fractional Laplacian that is defined as follows.
Let $\xO \subset \Ren$ be a bounded  domain, and   $\xl_i$ and  $\phi_i$ be the Dirichlet eigenvalues and
orthonormal  eigenfunctions of the Laplacian,
i.e. $-\Delta \phi_i = \xl_i \phi_i$ in $\xO$, with $\phi_i = 0$ on $\partial \xO$.
 Then,  for $f(x) = \sum c_i \phi_i(x)$ we  define
\be\la{1.frc}
(-\Delta)^{s} f(x) =  \sum_{i=1}^{\infty} c_i \xl_i^s  \phi_i(x),~~~~~~~~~~~~0<s<1 \ ,
\ee
in which case
\be\la{1.fr2}
((-\Delta)^{s} f , f )_{\xO} = \int_{\xO} f(x)~ (-\Delta)^{s} f(x) dx =  \sum_{i=1}^{\infty}  c_i^2 \xl_i^s.
\ee
In the sequel we will refer to this fractional Laplacian as the {\em spectral  fractional Laplacian}.
We then have

%\vspace{2mm}
\begin{theorem}\la{th13}({\bf Hardy \& Hardy-Sobolev-Maz'ya for  Spectral Fractional Laplacian}) \\
\noindent  Let $\frac12 \leq s <1$,  $n \geq 2$ and $\xO \subset \Ren$  be a bounded  domain. \\
\noindent  {\bf (i)} If in  addition  $\xO$ is
 such that
\be\la{c2b}
- \Delta d(x) \geq 0, ~~~~~~~~~~~~~ x \in \xO  \ ,
\ee
 then,  for all   $f \in C^{\infty}_{0}(\xO) $ there holds
\be\la{1.A22}
 ((-\Delta)^{s} f , f )_{\xO} \geq  d_{s}
 \int_{\xO} \frac{f^2(x)}{d^{2s}(x)}dx \ ,
\ee
with
\be\la{1.Ap2}
 d_{s} :=  \frac{2^{2s}  \Gamma^2 \left(\frac{3+2s}{4} \right)}
{\Gamma^2 \left(\frac{3-2s}{4} \right)} \ .
\ee
\noindent  {\bf(ii)}
 Suppose there exists a point $x_0 \in \partial \xO$ and $r>0$ such that the
part of the boundary $\partial \xO  \cap B(x_0,r)$ is $C^1$ regular. Then
\[
 d_{s} \geq
 \inf_{f \in C^{\infty}_{0}(\xO)}
 \frac{ ((-\Delta)^{s} f , f )_{\xO}}{\int_{\xO} \frac{f^2(x)}{d^{2s}(x)} dx}  \ .
\]
\noindent  {\bf(iii)} If $\xO$ is a  Lipschitz domain satisfying (\ref{c2b})
and $s \in (\frac12,1)$,
then there exists a positive constant $c$  such that
for all   $f \in C^{\infty}_{0}(\xO) $ there holds
\be\la{1.A22s}
 ((-\Delta)^{s} f , f )_{\xO} \geq  d_{s}
 \int_{\xO} \frac{f^2(x)}{d^{2s}(x)}dx + c
 \left(  \int_{\xO} | f(x)|^{\frac{2n}{n-2s}} dx \right)^{\frac{n-2s}{n}}  \ .
\ee

\end{theorem}

We next consider the fractional Laplacian associated to the  killed  stable processes  that appears in
 \cite{BD}, \cite{BBC}, \cite{SV1},
\cite{SV2}, \cite{PSV},  which from now on we will call it {\em Dirichlet fractional Laplacian}.
 The proper extension problem  involves test functions
 $u \in  C^{\infty}_{0}(\Ren \times \R)$ such  that
$u(x,0)=0$  in the complement of $\xO$, that is,  for $x \in \CC \xO$. For this fractional Laplacian, our assumption
on the domain $\xO$ is convexity instead of (\ref{c}). The reason for this is that our method requires subharmonicity
of the distance function in $\CC \xO$ which is equivalent to the convexity of $\xO$, see \cite{AK}.
 Our next result reads:

\begin{theorem}\la{th14} ({\bf Trace Hardy \& Trace Hardy-Sobolev-Maz'ya  II}) \\
\noindent Let $\frac12 \leq s <1$, $n \geq 2$  and $\xO \subsetneqq \Ren$  be a  domain.  \\
\noindent  {\bf (i)} If in addition  $\xO$ is convex
  then, for all   $u \in C^{\infty}_{0}(\Ren \times \R) $ such that $u(x,0)=0$  for $x \in \CC \xO$,  there holds
\be\la{1.B1}
 \dinr y^{1-2s} |\ana_{(x,y)} u(x,y)|^2 \dd  \geq \bar{k}_{s}  \int_{\xO} \frac{u^2(x,0)}{d^{2s}(x)}dx \ ,
\ee
with
\be\la{1.B2}
  \bar{k}_{s} :=  \frac{ 2^{1-2 s}  \Gamma^2(s +\frac12) \Gamma(1-s)}{\pi \Gamma(s)}  \ .
\ee

\noindent  {\bf(ii)}
 Suppose there exists a point $x_0 \in \partial \xO$ and $r>0$ such that the
part of the boundary $\partial \xO  \cap B(x_0,r)$ is $C^1$ regular. Then
\[
 \bar{k}_{s} \geq
\inf_{\scriptsize \begin{array}{c}
              u \in C^{\infty}_{0}(\Ren \times \R), \\
             u(x,0)=0,~ ~x \in \CC \xO
               \end{array}}
\frac{\dinr y^{1-2s} |\ana u|^2 \dd}{\int_{\xO} \frac{u^2(x,0)}{d^{2s}(x)} dx}  \ .
\]
% \inf_{u \in C^{\infty}_{0}(\Ren \times \R),~~ u(x,0)=0 ~x \in \CC \xO}
% \frac{\dinr y^{1-2s} |\ana u|^2 \dd}{\int_{\xO} \frac{u^2(x,0)}{d^{2s}(x)} dx}  \ .
In particular $ \bar{k}_{s}$ in (\ref{1.B1}) is the best constant. \\
\noindent  {\bf(iii)} If $\xO$ is a uniformly Lipschitz and  convex domain with finite inner radius
and $s \in (\frac12,1)$,
 then
 there exists a positive constant $c$, such that
 the following improvement holds true for all   $u \in C^{\infty}_{0}(\Ren \times \R) $
with  $u(x,0)=0$  for $x \in \CC \xO $:
\be\la{1.B3}
 \dinr y^{1-2s} |\ana_{(x,y)} u(x,y)|^2 \dd  \geq  \bar{k}_{s}
 \int_{\xO} \frac{u^2(x,0)}{d^{2s}(x)}dx + c
 \left(  \int_{\xO} | u(x,0)|^{\frac{2n}{n-2s}} dx \right)^{\frac{n-2s}{n}}  \ ,
\ee

\end{theorem}

Elementary manipulations show that
\[
\bar{d}_s = 2 \sin^2 \left(\frac{(2s+1) \pi}{4} \right) ~ \bar{k}_s    \ ,
\]
thus
\[
\bar{d}_s   > \bar{k}_s  \ ,   ~~~~~~{\rm for}~~~s \in (0,1) \ ,
\]
which implies in particular that the  best constants of Theorems \ref{th11} and \ref{th14} are different.

We next apply Theorem \ref{th14} to the  Dirichlet fractional Laplacian.
 In this  case,  for $f \in  C^{\infty}_{0}(\xO) $ we extend $f$ in all of $\Ren$ by setting
$f =0$ in $\CC \xO$ and use (\ref{1.2}).  In particular, the corresponding quadratic form  is
\bea
((-\Delta)^s_{D} f,f)_{\Ren} & = &
% \int_{\Ren} f~ (-\Delta)^{s} f dx =
\frac{c_{n,s}}{2} \int_{\R^n}\int_{\R^n}
\frac{|f(x)-f(\xi)|^2}{|x-\xi|^{n+2s}} dx d\xi \la{1.20}    \\
& = & \frac{c_{n,s}}{2} \left(\int_{\xO}\int_{\xO}
\frac{|f(x)-f(\xi)|^2}{|x-\xi|^{n+2s}} dx d\xi +
2  \int_{\xO}\int_{\CC \xO} \frac{|f(x)|^2}{|x-\xi|^{n+2s}} dx d\xi \right)  \ ,
\nonumber
\eea
with the constant $c_{n,s}$ as given by (\ref{1.cns}). We then have:

\begin{theorem}\la{th15}
 ({\bf  Hardy \& Hardy-Sobolev-Maz'ya  for  the Dirichlet  Fractional Laplacian}) \\
\noindent  Let $\frac12 \leq s <1$, $n \geq 2$  and $\xO \subsetneqq \Ren$  be a   domain.  \\
\noindent  {\bf (i)}If in addition $\xO$ is  convex,   then   for all   $f \in C^{\infty}_{0}(\xO) $ there holds
\be\la{1.B215}
 ((-\Delta)^{s}_{D} f , f )_{\Ren} \geq  \frac{\Gamma^2 \left(s+\frac12 \right)}{\pi}
 \int_{\xO} \frac{f^2(x)}{d^{2s}(x)}dx \ .
\ee
Equivalently, one has that
\bea\la{1.B315}
 \int_{\R^n}\int_{\R^n}
\frac{|f(x)-f(\xi)|^2}{|x-\xi|^{n+2s}} dx d\xi  \geq
 k_{n,s}
 \int_{\xO} \frac{f^2(x)}{d^{2s}(x)}dx     \ ,
\eea
where
\be\la{1.26}
k_{n,s} := \frac{2^{1-2s} \pi^{\frac{n-2}{2}} \Gamma(1-s) \Gamma^2(s + \frac12)}{s \Gamma( \frac{n+2s}{2})}  \ .
\ee
\noindent  {\bf(ii)}
 Suppose there exists a point $x_0 \in \partial \xO$ and $r>0$ such that the
part of the boundary $\partial \xO  \cap B(x_0,r)$ is $C^1$
regular. Then the Hardy constants  $ \frac{\Gamma^2
\left(s+\frac12 \right)}{\pi}$   in (\ref{1.B215})
  and $k_{n,s}$  in (\ref{1.B315}) are optimal.   \\
\noindent  {\bf(iii)} If  $\xO$ is a uniformly Lipschitz and
convex domain with finite inner radius and $s \in (\frac12,1)$,
  then there exists a positive constant $c$ such that
 for all   $f \in C^{\infty}_{0}(\xO) $ there holds
\be\la{1.B2s}
 ((-\Delta)^{s}_{D} f , f )_{\Ren} \geq  \frac{\Gamma^2 \left(s+\frac12 \right)}{\pi}
 \int_{\xO} \frac{f^2(x)}{d^{2s}(x)}dx + c
 \left(  \int_{\xO} | f(x)|^{\frac{2n}{n-2s}} dx \right)^{\frac{n-2s}{n}}  \ .
\ee
Equivalently, one has that
\bea\la{1.B3s}
 \int_{\R^n}\int_{\R^n}
\frac{|f(x)-f(\xi)|^2}{|x-\xi|^{n+2s}} dx d\xi  \geq
 k_{n,s}
 \int_{\xO} \frac{f^2(x)}{d^{2s}(x)}dx
+ c
 \left(  \int_{\xO} | f(x)|^{\frac{2n}{n-2s}} dx \right)^{\frac{n-2s}{n}}  \ .
\eea

\end{theorem}

The  case  where $\xO$ is the  half--space   $\xO = \Ren_{+} = \{ (x_1, \ldots, x_n): x_n >0 \}$  is of particular
interest see \cite{BD}, \cite{BBC}, \cite{FS}, \cite{D}, \cite{S}. In this case we obtain a stronger result that covers
the full range $s \in (0,1)$. More precisely we have:

\begin{theorem}\la{th16}
({\bf Half Space, Trace Hardy-Sobolev-Maz'ya  \& Fractional Hardy-Sobolev-Maz'ya II})  \\
 Let $0 < s <1$ and  $n \geq 2$. \\
\noindent
{\bf (i)}  Then
 for all   $u \in C^{\infty}_{0}(\Ren \times \R) $
with  $u(x,0)=0$, $x \in \Ren_{-}$, there holds
\be\la{1.C1}
 \dinr y^{1-2s} |\ana_{(x,y)} u(x,y)|^2 \dd  \geq  \bar{k}_{s}
 \int_{\Ren_{+}} \frac{u^2(x,0)}{x_n^{2s}}dx  \ ,
\ee
where
\[
  \bar{k}_{s} :=  \frac{ 2^{1-2 s}  \Gamma^2(s +\frac12) \Gamma(1-s)}{\pi \Gamma(s)}  \ ,
\]
is the best constant in (\ref{1.C1}). \\
\noindent  {\bf (ii)}   There exists a positive constant $c$, such that
for all   $u \in C^{\infty}_{0}(\Ren \times \R) $
with  $u(x,0)=0$, $x \in \Ren_{-}$, there holds
\be\la{1.C1s}
 \dinr y^{1-2s} |\ana_{(x,y)} u(x,y)|^2 \dd  \geq  \bar{k}_{s}
 \int_{\Ren_{+}} \frac{u^2(x,0)}{x_n^{2s}}dx + c
 \left(  \int_{\Ren_{+}} | u(x,0)|^{\frac{2n}{n-2s}} dx \right)^{\frac{n-2s}{n}}  \ ,
\ee
\noindent {\bf (iii)}
As a consequence, there exists a positive constant $c$  such  that
 for all   $f \in C^{\infty}_{0}(\Ren_{+}) $ there holds
\bea\la{1.C2}
 \int_{\R^n}\int_{\R^n}
\frac{|f(x)-f(\xi)|^2}{|x-\xi|^{n+2s}} dx d\xi  \geq
 k_{n,s}
 \int_{\Ren_{+}} \frac{f^2(x)}{x_n^{2s}}dx
+ c
 \left(  \int_{\Ren_{+}} | f(x)|^{\frac{2n}{n-2s}} dx \right)^{\frac{n-2s}{n}}  \ ,
\eea
where $k_{n,s}$ is given by (\ref{1.26}).

Or, equivalently,  for all   $f \in C^{\infty}_{0}(\Ren_{+}) $ there holds
\bea\la{1.C3}
 \int_{\R^n_{+}}\int_{\R^n_{+}}
\frac{|f(x)-f(\xi)|^2}{|x-\xi|^{n+2s}} dx d\xi  \geq  \xk_{n,s}
 \int_{\Ren_{+}} \frac{f^2(x)}{x_n^{2s}}dx
+ c
 \left(  \int_{\Ren_{+}} | f(x)|^{\frac{2n}{n-2s}} dx \right)^{\frac{n-2s}{n}}  \ ,
\eea
where
\[
\xk_{n,s} := \pi^{\frac{n-1}{2}} \frac{\Gamma(s+\frac12)}{s \Gamma(\frac{n+2s}{2})} \left[ \frac{2^{1-2s}}{\sqrt{\pi}}
\Gamma(1-s) \Gamma(s+\frac12)-1 \right] \ .
\]

%\frac{\pi^{\frac{n-1}{2}} \Gamma(s + \frac12)}{s  \Gamma(s + \frac{n}{2})} \
\end{theorem}

\noindent

 We note that the  Hardy--Sobolev--Maz'ya inequality (\ref{1.C2})
refers to the Dirichlet fractional Laplacian, associated to the killed stable processes whereas
inequality (\ref{1.C3}) is  associated to the censored stable processes.
The  Hardy constants $k_{n,s}$ and $\xk_{n,s}$ appearing in (\ref{1.C2}) and (\ref{1.C3}) respectively
 are optimal, as shown in \cite{BD}.
 The corresponding fractional Hardy inequality of (\ref{1.C3}) with
best constant,
 in the case of a convex domain $\xO$, that is,
\[
 \int_{\xO}\int_{\xO}
\frac{|f(x)-f(\xi)|^2}{|x-\xi|^{n+2s}} dx d\xi  \geq  \xk_{n,s}
 \int_{\xO} \frac{f^2(x)}{d(x)^{2s}}dx  , ~~~~~~~~~  f \in C^{\infty}_{0}(\xO)  \ ,
\]
has been established  for $s \in (\frac12,1)$  in \cite{LS}.
 The question of obtaining   a Hardy--Sobolev--Maz'ya inequality for the
half space was raised in \cite{FS} and was answered positively in \cite{S}, \cite{D},
 but only for the range $s \in (\frac12,1)$.

For other type of trace Hardy inequalities we refer to \cite{DDM}  and \cite{AFV}.
We finally note that fractional Sobolev inequalities play an important role in many other directions, see e.g.,
\cite{BBM},  \cite{CG}, \cite{MS}, \cite{N}.

\setcounter{equation}{0}
\section{The  Trace Hardy inequality I}\la{sec2}
In this section we will prove the trace Hardy inequality contained in Theorem \ref{th11}. We first recall
 the definition of
a uniformly Lipschitz domain $\xO$; see section 12 of \cite{L}. We  note that Stein calls such a domain
 minimally smooth, see section 3.3 of \cite{St}.

A domain $\xO$ is called uniformly Lipschitz if there exist $\xe>0$, $L>0$, and   $M  \in \N$ and a locally
 finite countable  cover $\{U_i\}$ of $\partial \xO$ with the following properties: \\
\noindent (i) If $x \in \partial \xO$ then $B(x,\xe) \subset U_i$ for some $i$. \\
\noindent (ii) Every point of $\R^n$ is contained in at  most $M$ $U_i$'s.  \\
\noindent (iii) For each $i$ there exist  local coordinates $y=(y',y_n) \in \R^{n-1} \times \R$ and
a Lipschitz function $f: \R^{n-1} \ra \R$, with  $Lipf \leq  L$ such that
\[
U_i \cap \xO= U_i \cap \{(y',y_n) \in \R^{n-1} \times \R: y_n > f(y') \}.
\]
Under the uniformly Lipschitz  assumption on $\xO$  the extension operator is defined in $W^{1,p}(\xO)$,
for all $p \geq 1$.  We also note that when $\xO$ is a bounded domain  the above definition reduces to
$\xO$ being Lipschitz.

In the sequel we set  $a=1-2s$.  Since  $ 0< s <1$ we also have   $-1 <  a < 1$.
 We first establish
the following useful identity:

\begin{lemma}\la{lem21}
Suppose that $a \in (-1,1)$ and let   $u \in C^{\infty}_{0}(\xO \times \R)$ and
  $\phi \in  C^{2}(\xO \times (0, \infty)) \cap
 C(\bar{\xO} \times [0, \infty))$  is such that  $\phi(x,y)>0$ in  $\xO \times [0, \infty)$,
$\phi(x,y)=0$ in  $\partial \xO \times (0, \infty)$,
\[
| y^a \frac{\phi_{y}(x,y) }{\phi(x,y)}| \leq V(x),~~~y \in (0,1), ~~~x \in \xO,~~~~~~~~ 0 \leq V(x) \in L^{1}_{loc}(\xO),
\]
and for  a.e.  $x \in \xO$,  the following limit exists:
\[
\lim_{y \ra 0^+} \left( y^a \frac{\phi_{y}(x,y) }{\phi(x,y)} \right) \ .
\]
We also require that  the following integrals are finite
\[
 \dino y^a  \frac{|\ana \phi|^2}{\phi^2} u^2 \dd, ~~~~~~~~~~
 \dino  \frac{|{\rm div}(y^a \ana \phi)|}{\phi} u^2 \dd  \ .
\]
We then have the identity:
\bea\la{iden}
 \dino y^a |\ana u|^2 \dd     =  & -  &
  \int_{\xO}\lim_{y \ra 0^+} \left( y^a \frac{\phi_{y} }{\phi} \right) u^2(x,0) dx
+ \dino y^a |\ana u - \frac{\ana \phi}{\phi} u|^2 \dd  \nonumber \\
&  - &
 \dino \frac{{\rm div}(y^a \ana \phi)}{\phi} u^2 \dd.
\eea
\end{lemma}

{\em Proof:} Expanding the square and integrating by parts we compute for $\xe>0$,
%\subsection{The identity}\la{sec2.1}
\bean
\int_{\xe}^{\infty}\int_{\xO} y^a |\ana u - \frac{\ana \phi}{\phi} u|^2 \dd  &  = &
 \int_{\xe}^{\infty}\int_{\xO}  y^a  \left( |\ana u|^2 +
\frac{|\ana \phi|^2}{\phi^2} u^2 - \frac{\ana \phi}{\phi} \ana u^2 \right) \dd
  \\
 =  \int_{\xe}^{\infty}\int_{\xO} y^a |\ana u|^2 \dd  &  +  &
  \int_{\xe}^{\infty}\int_{\xO} \frac{{\rm div}(y^a \ana \phi)}{\phi} u^2 \dd + \int_{\xO}
  \xe^a \frac{\phi_{y}(x,\xe)}{\phi(x,\xe)}  u^2(x,\xe) dx \ .
\eean
We then pass to limit $\xe \ra 0$ and the  result follows easily.

\finedim

We will use Lemma \ref{lem21} with the following choice:
 $\phi(x,y)= d^{-\frac{a}{2}}(x) A \left( \frac{y}{d(x)} \right)$ for
 $y >0$, $x \in \Omega$.  The function  $A$ solves the following boundary value problem
\be\la{ode1}
(t^3+t)A'' + (a + t^2(2+a))A' +  \frac{(2+a)a}{4} t A  =  0, ~~~~t>0,
\ee
with
\be\la{2bc}
A(0)=1, ~~~~~~~~~~~~~ \lim_{t \rft} A(t) = 0 \ .
\ee
Equation (\ref{ode1}) can also  be written in divergence form as
\be\la{ode1div}
(t^a(1+t^2)A')' + \frac{(2+a)a}{4} t^a A = 0.
\ee

From now on we will use the following notation:
\[
f \sim g, ~~~~~~~~in ~~~U,
\]
whenever there exist positive constants $c_1$, $c_2$, such that
\[
c_1 g \leq f \leq c_2 g, ~~~~~~~~~ in ~~U  \ .
\]

We then have the following

\begin{prop}\la{prop22}
Suppose that $a \in (-1,1)$.
The boundary value problem (\ref{ode1}), (\ref{2bc}) has a positive decreasing  solution $A$ with the
following properties:  \\
\noindent (i) There exists a positive constant $\bar{d}_s$ such that
\[
      \lim_{t \ra 0^{+}} t^a A'(t) = - \bar{d}_s \ ,
\]
with
\[
\bar{d}_s = \frac{(1-a) \Gamma \left(\frac{1+a}{2} \right) \Gamma^2 \left(\frac{4-a}{4} \right)}
{\Gamma^2 \left(\frac{2+a}{4} \right) \Gamma \left(\frac{3-a}{2} \right)}
=  \frac{2s \Gamma \left(1-s \right) \Gamma^2 \left(\frac{3+2s}{4} \right)}
{\Gamma^2 \left(\frac{3-2s}{4} \right) \Gamma \left( 1+s \right)}.
\]
(ii) For all $t>0$,
\bean
   A(t)    &  \sim  & (1+t^2)^{-\frac{2+a}{4}}   \ ,   \\
A^{\prime}(t)   &  \sim  &  - t^{-a} (1+t^2)^{-\frac{4-a}{4}}  \ .
\eean
Moreover,
\[
\lim_{t \rft} \frac{tA^{\prime}(t)}{A(t)} = -\frac{2+a}{2}  \ .
\]
(iii) There holds:
\be\la{dsben}
\bar{d}_s = \int_{0}^{\infty} t^{a}(1+t^2) (A')^2 dt -  \frac{(2+a)a}{4}  \int_{0}^{\infty} t^a A^2 dt,
\ee
(iv) In case $a \in (-1,0]$, we  have
\[
 t A'(t) + \frac{a}{2} A(t)  \leq 0 \ .
\]
Moreover  for  $a \in (-1,0)$ and all $t>0$ we have
\[
 t A'(t) + \frac{a}{2} A(t)   \sim - A(t) \ .
\]

\end{prop}

\noindent  {\em Proof:}
We change variables in (\ref{ode1})  by $z=-t^2$ and define $B(z)$ such that $A(t)=B(-t^2)$, whence $A_t = -2t B_z$ and
$A_{tt}= -2B_z + 4t^2 B_{zz}$. It then follows that $B(z)$ satisfies the Gauss hypergeometric equation
\[
z(1-z) B '' + \left( \frac{1+a}{2} -\frac{3+a}{2}z  \right) B' - \frac{a(2+a)}{16}B = 0,~~~~~~~~~-\infty < z <0,
\]
whose general solution is given by
\[
B(z) = C_1 F_1 \left( \frac{a}{4}, \frac{2+a}{4}, \frac{1+a}{2}; z\right) +
 C_2 z^{\frac{1-a}{2}} F_2 \left( \frac{2-a}{4}, \frac{4-a}{4}, \frac{3-a}{2}; z \right) \ ;
\]
see \cite{AS}, Section 15.5  as well as 15.1 for  the definition and basic properties of the function $F$.
 It follows that
\be\la{sola}
A(t) =C_1 F_1 \left( \frac{a}{4}, \frac{2+a}{4}, \frac{1+a}{2}; -t^2 \right) +
 C_2 t^{1-a}e^{\frac{i \pi (1-a)}{2}} F_2 \left( \frac{2-a}{4}, \frac{4-a}{4}, \frac{3-a}{2}; -t^2 \right).
\ee
Since $F(\xa, \xb, \xg ; 0)=1$ for any $\xa$, $\xb$, $\xg$, the condition $A(0)=1$ implies that $C_1=1$.
We then have
\bea\la{dsb2}
\bar{d}_s &  = &  - \lim_{t \ra 0^+} t^a A'(t) \nonumber \\
 &  = &  - \lim_{t \ra 0^+}t^a   ( -2t F_1' + (1-a) C_2 e^{\frac{i \pi (1-a)}{2}} t^{-a}F_2
-2 C_2  t^{2-a}e^{\frac{i \pi (1-a)}{2}} F_2' ) \nonumber \\
&  = &  -(1-a) C_2 e^{\frac{i \pi (1-a)}{2}}.
\eea
In the above calculation we have also used the fact that
\[
F'(\xa, \xb, \xg ; z) = \frac{d}{dz}F(\xa, \xb, \xg ; z) = \frac{\xa \xb}{\xg} F(\xa+1, \xb+1, \xg+1 ; z).
\]
We next compute the behavior of $A$ at infinity. To this end we will  use the inversion formula, valid for
any $\xa$, $\xb$, $\xg$ and $|arg(-z)|< \pi$:
\bean
F(\xa, \xb, \xg ; z) = & &
 \frac{\xG(\xg) \xG(\xb-\xa)}{\xG(\xb) \xG(\xg-\xa)}(-z)^{-\xa} F\left( \xa,~ 1-\xg+\xa,~ 1-\xb+\xa; \frac{1}{z} \right)  \\
& + &
\frac{\xG(\xg) \xG(\xa-\xb)}{\xG(\xa) \xG(\xg-\xb)}(-z)^{-\xb} F\left( \xb,~ 1-\xg+\xb,~ 1-\xa+\xb; \frac{1}{z} \right).
\eean
We then calculate
\bean
 \lim_{t \rft} t^{\frac{a}{2}} A(t) = \frac{ \Gamma \left(\frac{1+a}{2} \right)
\Gamma \left(\frac{1}{2} \right)} {\Gamma^2 \left(\frac{2+a}{4} \right)} + C_2  e^{\frac{i \pi (1-a)}{2}}
 \frac{ \Gamma \left(\frac{3-a}{2} \right)\Gamma \left(\frac{1}{2} \right)}{\Gamma^2 \left(\frac{4-a}{4} \right)}.
\eean
To make this limit equal to zero  we choose
\[
C_2 = -  e^{-\frac{i \pi (1-a)}{2}} \frac{ \Gamma \left(\frac{1+a}{2} \right)
\Gamma^2 \left(\frac{4-a}{4} \right)} {\Gamma^2 \left(\frac{2+a}{4} \right)\Gamma \left(\frac{3-a}{2} \right)}.
\]
Combining this with (\ref{dsb2}) we conclude
\be\la{ctrval}
\bar{d}_s = \frac{(1-a) \Gamma \left(\frac{1+a}{2} \right) \Gamma^2 \left(\frac{4-a}{4} \right)}
{\Gamma^2 \left(\frac{2+a}{4} \right) \Gamma \left(\frac{3-a}{2} \right)}
=  \frac{2s \Gamma \left(1-s \right) \Gamma^2 \left(\frac{3+2s}{4} \right)}
{\Gamma^2 \left(\frac{3-2s}{4} \right) \Gamma \left( 1+s \right)}.
\ee

At this point both constants $C_1$, $C_2$, in (\ref{sola}) have been identified. After some lengthy but
 straightforward calculations
we find that as  $t \rft$
\be\la{asa}
A(t) \sim t^{-\frac{2+a}{2}},~~~~~~~~~~ A'(t) \sim t^{-\frac{4+a}{2}} \ .
\ee
In addition we  get
\[
\lim_{t \rft} \frac{tA^{\prime}(t)}{A(t)} = -\frac{2+a}{2}  \ .
\]
Using (\ref{ode1div}) and the above asymptotics, we easily conclude that the solution $A$ is energetic, that is,
\[
 \int_{0}^{\infty} t^{a}(1+t^2) (A')^2 dt  +   \int_{0}^{\infty} t^a A^2 dt  < \infty \ .
\]
Multiplying (\ref{ode1div})  by A and integrating by parts in $(0,\infty)$ we arrive at (\ref{dsben})

To prove the positivity and monotonicity of $A$ we next change variables by:
\[
B(s) = (1+t^2)^{\frac{a}{4}} A(t), ~~~~~~~~~~~s=1/t \ .
\]
It follows that $B$ satisfies the equation
\[
(1+s^2)^2 B''+
 (2-a)s(1+s^2) B' -\frac{a^2}{4} B =0,  ~~~~~~~~ s \in (0, +\infty)  \ ,
\]
with $B(0)=0$ and $B(+\infty)=1$.
 A standard maximum principle argument shows that $B$ is positive. Consequently $A$
is positive and the monotonicity of $A$ follows easily.

The positivity and monotonicity  of $A$ in connection with the asymptotics of $A$ yield easily
part (ii) of the Proposition.

Part (iv) follows easily from the monotonicity of $A$ and part (ii).

\finedim

Using the asymptotics of $A(t)$, from the previous  Proposition we easily obtain  the following uniform
asymptotics for $\phi$
\begin{lemma}\la{lemasf}
 Suppose  $a \in (-1,1)$  and let $\phi$ be given by
\[
\phi(x,y)= d^{-\frac{a}{2}}(x) A \left( \frac{y}{d(x)} \right), ~~~~~~~~
 y >0,~~~ x \in \Omega \subset \R^n  \ ,
\]
where $A$ solves (\ref{ode1}), (\ref{2bc}).

\noindent ({\bf i})
 Then
\[
\phi(x,y) \sim \frac{d}{(d^2+y^2)^{\frac{2+a}{4}}}  \ ,  ~~~~ y >0,~~~  x \in \Omega \ .
\]
Concerning the gradient of $\phi$,  for $a \in (-1,0]$ we have
\[
|\ana_{(x,y)} \phi(x,y)| \sim \frac{1}{(d^2+y^2)^{\frac{2+a}{4}}} \ ,  ~~~~ y >0,~~~  x \in \Omega \ ,
\]
whereas for  $a \in (0,1)$
\[
|\ana_{(x,y)} \phi(x,y)| \sim \frac{y^{-a}}{(d^2+y^2)^{\frac{2-a}{4}}} \ ,  ~~~~ y >0,~~~  x \in \Omega \ .
\]

\noindent ({\bf ii})
If     $\xO$ satisfies  $-\Delta d(x) \geq 0$ for $x \in \xO$, then   for    $a \in (-1,0)$
\[
-{\rm div}(y^a \ana \phi) \phi  \sim  \frac{y^{a}}{(d^2+y^2)^{\frac{2+a}{2}}} (-d \Delta d) \ ,~~~~ y >0,~~~  x \in \Omega \ ,
\]
whereas for $a=0$,
\[
-{\rm div}( \ana \phi) \phi  \sim  \frac{y}{(d^2+y^2)^{\frac{3}{2}}} (-d \Delta d) \ ,~~~~ y >0,~~~  x \in \Omega \ .
\]

\end{lemma}

We are now ready to give the proof of Theorem \ref{th11}.

\noindent
{\em Proof of Theorem \ref{th11}  part (i) and (ii):} We assume that $s \in [\frac12,1)$ or equivalently $a \in (-1,0]$.
 We will use Lemma  \ref{lem21}
with   the test function $\phi$ given by
\[
\phi(x,y)= d^{-\frac{a}{2}}(x) A \left( \frac{y}{d(x)} \right), ~~~~~~~~
 y >0,~~~ x \in \Omega \subset \R^n  \ ,
\]
where $A$ solves (\ref{ode1}), (\ref{2bc}).  Using Proposition \ref{prop22}   and Lemma  \ref{lemasf}  we see
that all hypotheses of Lemma \ref{lem21} are satisfied. In particular,
 for $t = \frac{y}{d}$ we compute, for $x \in \xO$,
\bea\la{trco}
-\lim_{y \ra 0^+} \left( y^a \frac{\phi_{y} }{\phi} \right) &  =  &
-\lim_{y \ra 0^+} \left(  t^a \frac{A'(t)}{d^{1-a} A(t)} \right) =
\frac{1}{d^{1-a}(x)}  \lim_{t \ra 0^+} \left(-  \frac{ t^a A'(t)}{A(t)} \right)   \nonumber \\
&  =  & \frac{ \bar{d}_s }{d^{1-a}(x)}  \ .
\eea
We also have
\bean
- {\rm div}(y^a \ana \phi)&  = & - y^{a-1} d^{-1-\frac{a}{2}}\left[
(t^3+t)A'' + (a + t^2(2+a))A' +  \frac{(2+a)a}{4} t A  \right]   \\
&  &   - y^{a-1} d^{-1-\frac{a}{2}}\left[ (-d \Delta d) \left(t^2 A' + \frac{at}{2} A \right) \right]  \\
&  = & - y^{a-1} d^{-1-\frac{a}{2}} \left[(-d \Delta d) \left(t^2 A' + \frac{at}{2} A \right) \right] \ ,
\eean
therefore,
\[
- {\rm div}(y^a \ana \phi)  \geq 0 \ , ~~~~~x \in \xO, ~~~~y>0.
\]
From Lemma \ref{lem21} we get \bea\la{id2}
 \dino y^a |\ana u|^2 \dd     \geq  &   &
\bar{d}_s  \int_{\xO}  \frac{ u^2(x,0)}{d^{1-a}(x)} dx
+ \dino y^a |\ana u - \frac{\ana \phi}{\phi} u|^2 \dd  \nonumber \\
&  - &
 \dino \frac{{\rm div}(y^a \ana \phi)}{\phi} u^2 \dd  \ ,
\eea
from which the trace Hardy inequality follows directly.  This relation will be used later on,
 in Sections \ref{sec4}  and \ref{sec5}
 to obtain the Sobolev term as well.

We continue with the proof of the optimality of
 the Hardy constant  $\bar{d}_s$.
Let
\be\la{rq}
Q[u] := \frac{\dino y^a |\ana u|^2 \dd}{\int_{\xO} \frac{u^2(x,0)}{d^{1-a}(x)} dx}=:\frac{N[u]}{D[u]}.
\ee
We have that $Q[u] \geq \bar{d}_s$. Here we will show that there exists a sequence of functions $u_{\xe}$ such that
$\lim_{\xe \ra 0}Q[u_{\xe}] =  \bar{d}_s$, and therefore $\bar{d}_s$ is the best constant.

We first assume for simplicity that the boundary of $\xO$  is flat in a neighborhood $V$  of a point $x_0 \in \partial \xO$.
 The  neighborhood of the point $x_0$ is assumed to  contain a ball centered
at $x_0$ with radius, say,  $3\xd$.
Locally around $x_0$ the boundary is given by $x_n=0$, whereas the interior of $\xO$ corresponds to $x_n>0$. We also
write $x=(x', x_n)$. Clearly, for $x \in \xO \cap V$  we have that $d(x) = x_n$.

We next define two suitable cutoff functions.
 Let $\psi(x') \in C_0^{\infty}(B_{\xd})$,
 where $B_{\xd} \subset \partial \xO \subset \R^{n-1}$ is the ball centered at $x_0$ with radius
$\xd$. Also the nonnegative function
 $h(x_n) \in C^{\infty}(\R^+)$ is such that $h(x_n) =0$  for $x_n \geq 2 \xd$ and $h(x_n) =1$ for $0 \leq x_n \leq \xd$.
We will use the following test function:
\be\la{ue}
u_{\xe}(x',x_n,y)   = \left\{ \begin{array}{ll}
 h(x_n) \psi(x')x_n^{-\frac{a}{2}} A(\frac{y}{x_n}),    &  ~~~~~~~~~  y \geq \xe  \\
  h(x_n) \psi(x')x_n^{-\frac{a}{2}} A(\frac{\xe}{x_n}), &~~~~~0\leq y< \xe.\\
\end{array} \right.
\ee We have that \be\la{qe} Q[u_{\xe}] = \frac{ \int_0^{+\infty}
dy \int_0^{2 \xd} dx_n \int_{B_{\xd}} dx' y^a |\ana u_{\xe}|^2} {
\int_0^{2 \xd} dx_n  \int_{B_{\xd}}  dx'
\frac{u_{\xe}^2}{x_n^{1-a}}} = \frac{N[u_{\xe}]}{D[u_{\xe}]}. \ee
Concerning the denominator we compute \bea\la{deno} D[u_{\xe}] &
= &  \int_{B_{\xd}}\psi^2(x')  dx' ~ \int_0^{ \xd} x_n^{-1}
A^2(\frac{\xe}{x_n}) dx_n +O_{\xe}(1)
\nonumber \\
 &  = &  \int_{B_{\xd}}\psi^2(x')  dx' ~ \int_{\xe/\xd}^{+\infty} \frac{A^2(t)}{t} dt +O_{\xe}(1).
\eea

We next calculate the numerator. At first we break $N$ into two pieces:
\[
N[u_{\xe}] = \int_{0}^{\xe} dy + \int_{\xe}^{+\infty} dy =:N_1[u_{\xe}]+ N_2[u_{\xe}].
\]
Using the specific form of $u_{\xe}$ and elementary estimates  we
calculate: \bean N_2[u_{\xe}]&  = &  \int_{B_{\xd}}  \psi^2(x')
dx'    \int_{\xe}^{+\infty} dy \int_0^{ \xd} dx_n  ~
 \frac{y^a}{x^{a+2}_n}
  \left[
\left( -\frac{a}{2} A(\frac{y}{x_n})- \frac{y}{x_n} A'(\frac{y}{x_n}) \right)^2
+    A^{'2}(\frac{y}{x_n}) \right]
 \\
& + &  \int_{B_{\xd}}  |\ana \psi(x')|^2 dx'
\int_{\xe}^{+\infty} dy \int_0^{ \xd} dx_n  ~ y^a  x_n^{-a}
 A^2(\frac{y}{x_n})   + O_{\xe}(1)  \\
&  =: & N_{21}[u_{\xe}] +N_{22}[u_{\xe}] + O_{\xe}(1).
\eean
We note that as $\xe \ra 0$,
\bean
N_{22}[u_{\xe}]
 &  =  &
 \int_{B_{\xd}}  |\ana \psi(x')|^2 dx'
 \int_0^{ \xd} x_n   \int_{\xe/x_n}^{+\infty} t^a A^2(t) dt d x_n
  \\
%& =  &   \frac12    \xd^2   \int_{B_{\xd}}  |\ana \psi(x')|^2 dx'  \int_{0}^{+\infty} t^a A^2(t) dt  \\
& =  & O_{\xe}(1). \eean Concerning $N_{21}[u_{\xe}]$, changing
variables by   $ t = \frac{y}{x_n}$   we   write: \bean
N_{21}[u_{\xe}] &  = &  \int_{B_{\xd}}  \psi^2(x') dx'
\int_{\xe}^{+\infty} \frac{dy}{y}
 \int_{y/\xd}^{+\infty} \left[ t^a A^{'2}(t) + t^a \left(\frac{a}{2}A(t)+ tA'(t)\right)^2 \right] dt \\
 &  = &  \int_{B_{\xd}}  \psi^2(x') dx'   \int_{\xe}^{+\infty} \frac{dy}{y}  \int_{y/\xd}^{+\infty} \left[
t^a(1+t^2) A^{'2} + at^{1+a}AA' + \frac{a^2}{4}t^a A^2 \right]dt.
\eean
Integrating by parts the term containing the factors $AA'$ and then using the
 equation satisfied by $A$ (cf (\ref{ode1div})) we get
\bean
\int_{y/\xd}^{+\infty}  \hspace{-8mm}  & &   \left[
t^a(1+t^2) A^{'2}   +   at^{1+a}AA' + \frac{a^2}{4}t^a A^2 \right]dt \hspace{4cm}   \\
&  = &  \int_{y/\xd}^{+\infty}\left[
t^a(1+t^2) A^{'2}  - \frac{a(2+a)}{4}t^a A^2 \right]dt
+ \frac12 at^{1+a}A^2(t)\Large|_{t=\frac{y}{\xd}} \\
&  = & -t^a(1+t^2)A(t) A'(t) \Large|_{t=\frac{y}{\xd}} + \frac12 at^{1+a}A^2(t)\Large|_{t=\frac{y}{\xd}},
\eean
whence,
\[
N_{21}[u_{\xe}] =  - \int_{B_{\xd}}  \psi^2(x') dx'  \int_{\xe/\xd}^{+\infty}
\frac{1}{t} t^a(1+t^2)A(t)A'(t) dt +  O_{\xe}(1).
\]
It is not difficult to show that $N_1[u_{\xe}] = O_{\xe}(1)$, and therefore $N[u_{\xe}] = N_{21}[u_{\xe}] +  O_{\xe}(1)$.
Using also (\ref{deno}) we can form the quotient
\bea\la{telos}
\lim_{\xe \ra 0}   Q[u_{\xe}] &  =  & \lim_{\xe \ra 0}
 \frac{ - \int_{B_{\xd}}  \psi^2(x') dx'  \int_{\xe/\xd}^{+\infty}
\frac{1}{t} t^a(1+t^2)A(t)A'(t) dt +  O_{\xe}(1)}
{\int_{B_{\xd}}  \psi^2(x') dx' \int_{\xe/\xd}^{+\infty} \frac{A^2(t)}{t} dt +O_{\xe}(1)} \nonumber    \\
& = &  \lim_{\xe \ra 0} \frac{ -   \int_{\xe/\xd}^{+\infty}
\frac{1}{t} t^a(1+t^2)A(t)A'(t) dt}{ \int_{\xe/\xd}^{+\infty} \frac{A^2(t)}{t} dt}  \nonumber    \\
& = & - \lim_{\xs \ra 0}  \frac{ \xs^a (1+\xs^2) A'(\xs)}{A(\xs)}  \nonumber      \\
& = &  \bar{d}_s,
\eea
where we  used L'Hopital's rule and then part (i) of Proposition \ref{prop22}.

Let us now consider the general case.  We assume that $\partial \xO$ is $C^1$ in a neighborhood of a point  $\bx_0$,
which we take to be the origin
  $0  \in \partial \xO$. Thus locally  $\partial \xO$, is the graph of a function $\bx_n= \xg(\bx')$, with $\xg(0)=0$
and $\ana \xg(0)=0$. We also assume that the interior of $\xO$ corresponds to $\bx_n > \xg(\bx')$. Then the following
change of coordinates straightens the boundary in a neighborhood of the origin:
 $x_i = \bx_i$, $i=1,2,\ldots,n-1$, and $x_n = \bx_n -\xg(\bx')$; see e.g. \cite{E}, Appendix C.
 We assume that inside  the ball $B(0, 3 \xd)$
 (in the $x$-space) the image of $\partial \xO$ is flat. We then consider the test function
 $v_{\xe}(\bx,y)= u_{\xe}(x,y)$. Clearly  $v_{\xe}(\bx,y)$ is zero away from  a neighborhood of the origin, say $U$, and
elementary calculations show that
\[
\ana_{\bx}v_{\xe} = \ana_{x}u_{\xe} - u_{\xe,x_n} \ana_{\bx} \xg(\bx'),
\]
whence,
\[
|\ana_{\bx}v_{\xe} - \ana_{x}u_{\xe} | \leq   | \ana_{\bx} \xg(\bx')| |\ana_{x}u_{\xe}| =o_{\xd}(1) |\ana_{x}u_{\xe}|.
\]
It then follows that
\[
|\ana_{\bx}v_{\xe}| = |\ana_{x}u_{\xe}| (1 + o_{\xd}(1)).
\]
On the other hand, for $\bx  \in U$ and $d(\bx)={\rm dist}(\bx,\partial \xO)$,  we have that
\[
d(\bx)  =  (\bx_n-\xg(\bx')) (1+ | \ana_{\bx} \xg(\bx')|^2)^{1/2}  =   x_n (1 + o_{\xd}(1)).
\]
We finally note that the Jacobian of the above transformation is one and therefore $dx = d \bx$. We then compute
\[
Q[v_{\xe}(\bx,y)] = Q[u_{\xe}(x,y)] (1 + o_{\xd}(1)),
\]
where $Q[u_{\xe}(x,y)]$ is given in (\ref{qe}).
Since $\xd$ can be taken as small as we like the result follows easily, using the calculations from the flat case.

\finedim

\setcounter{equation}{0}
\section{The  Trace Hardy inequality  II}\la{sec3}

In this section we will prove the trace Hardy inequality contained in Theorem \ref{th14}.
We first establish the analogue of Lemma  \ref{lem21}:

\begin{lemma}\la{lem31}
Suppose that $a \in (-1,1)$ and let   $u \in C^{\infty}_{0}(\Ren \times \R)$ such that  $u(\cdot, 0) \in C^{\infty}_{0}(\xO)$.
Let  $\phi \in  C^{2}(\Ren \times (0, \infty)) \cap
 C(\Ren \times [0, \infty))$  is such that  $\phi(x,y)>0$ in  $\Ren \times [0, \infty)$,
$\phi(x,0)=0$ in  $x \in \CC \xO$,
\[
| y^a \frac{\phi_{y}(x,y) }{\phi(x,y)}| \leq V(x),~~~y \in (0,1), ~~~x \in \Ren,~~~~~~~~ 0 \leq V(x) \in L^{1}_{loc}(\Ren) \ .
\]
Moreover  for  a.e.  $x \in \xO$,  the following limit exists:
\[
\lim_{y \ra 0^+} \left( y^a \frac{\phi_{y}(x,y) }{\phi(x,y)} \right) \ .
\]
We also require that  the following integrals are finite
\[
 \dinr y^a  \frac{|\ana \phi|^2}{\phi^2} u^2 \dd, ~~~~~~~~~~
 \dinr  \frac{|{\rm div}(y^a \ana \phi)|}{\phi} u^2 \dd  \ .
\]
We then have the identity:
\bea\la{iden2}
 \dinr y^a |\ana u|^2 \dd     =  & -  &
  \int_{\xO}\lim_{y \ra 0^+} \left( y^a \frac{\phi_{y} }{\phi} \right) u^2(x,0) dx
+ \dinr y^a |\ana u - \frac{\ana \phi}{\phi} u|^2 \dd  \nonumber \\
&  - &
 \dinr \frac{{\rm div}(y^a \ana \phi)}{\phi} u^2 \dd.
\eea
\end{lemma}

The proof of this Lemma is quite similar to the proof of Lemma \ref{lem21} and we omit it.

This time we will choose   the test function to be  of the form
\be\la{fi2}
\phi(x,y)   = \left\{ \begin{array}{ll}
 (y^2+d^2)^{-\frac{a}{4}} B(\frac{d}{y}),  &  ~~~~~~~  x \in \xO,~~~ y >0     \\
  (y^2+d^2)^{-\frac{a}{4}} B(-\frac{d}{y}), & ~~~~~~~ x \in \CC \xO,~~ y >0    \\
\end{array} \right.
\ee
where function
 $B$ is  the solution of the following  boundary value problem
\be\la{ode4}
(1+t^2)^2 B''+
 (2-a)t(1+t^2) B' -\frac{a^2}{4} B =0,  ~~~~~~~~ t \in (-\infty, +\infty)  \ ,
\ee
complemented with the conditions
\be\la{bc2}
B(-\infty) = 0, ~~~~~~~~~B(+\infty) = 1.
\ee
We note that this can be written in divergence form as
\be\la{odediv2}
((1+t^2)^{1-\frac{a}{2}} B'(t))' - \frac{a^2}{4}(1+t^2)^{-1-\frac{a}{2}} B(t) =0, ~~~~~t \in \R.
\ee
We next collect some properties of $B$ that will be used later on.

\begin{prop}\la{prop32}
Suppose that $a \in (-1,1)$.
The boundary value problem (\ref{ode4}), (\ref{bc2}) has a positive increasing  solution $B$ with the
following properties:  \\
\noindent (i) There exists a positive constant $\bar{k}_s$ such that
\be\la{tctr}
 \lim_{t \rft} (1+t^2)^{\frac{2-a}{2}} B'(t) =:  \tC  \ ,
\ee
where
\[
\tC
= \frac{2^a \Gamma^2(\frac{2-a}{2}) \Gamma(\frac{1+a}{2})}{  \pi \Gamma(\frac{1-a}{2})}
=\frac{2^{1-2 s}}{\pi} \frac{\Gamma^2(s +\frac12) \Gamma(1-s)}{\Gamma(s)} \ .
\]

\noindent (ii) We have
\bean
B(t) &  \sim  &  1,~~~~~~~~~~~~~~~ t>0  \\
B(t) &  \sim  &  (1+t^2)^{-\frac{1-a}{2}} ~~~~~~~~~ t<0 \ ,  \\
B^{\prime}(t) & \sim & (1+t^2)^{-\frac{2-a}{2}} ~~~~~~~~~ t \in \R \ .
\eean
(iii) There holds:
\[
 \tC =  \int_{-\infty}^{+\infty}\left[ (1+t^2)^{1-\frac{a}{2}} B^{'2}(t) +
 \frac{a^2}{4}
(1+t^2)^{-1-\frac{a}{2}} B^{2}(t) \right] dt \ .
 \]
(iv) In case $a \in (-1,0]$, we  have
\[
(1+t^2)B^{\prime}(t) - \frac{a}{2}t B(t) >0, ~~~~~~~~~t \in \R \ .
\]
Moreover for  $a \in (-1,0)$
\[
(1+t^2)B^{\prime}(t) - \frac{a}{2}t B(t) \sim (1+t^2)^{\frac{1}{2}}, ~~~~~~~~~ t >0 \ .
\]

\end{prop}

\noindent {\em Proof:} When $a=0$ the ODE can be easily  solved by a straightforward integration. For
the general case  we first change variables by $B(t) = (1+t^2)^{\frac{a}{4}} f(t)$ to obtain
\[
(1+t^2)f''+ 2t f'+ \frac{a(2-a)}{4} f = 0.
\]
We next change variables by $g(z)=f(t)$,  $z=it$, so that $g$ satisfies the equation
\be\la{ast}
(1-z^2) g''-2zg' + \nu(\nu+1)g = 0,     ~~~~~~~~~~\nu=-\frac{a}{2}.
\ee
The solution of this is given  in \cite{AS}, Section 8.1:
\be\la{solg}
g(z) = \left\{ \begin{array}{ll}
 C_1^{+} P_{\nu}(z) + C_2^{+} Q_{\nu}(z),  &  ~~~~~~~ {\rm Im} z >0,  \\
 C_1^{-} P_{\nu}(z) + C_2^{-} Q_{\nu}(z),  &  ~~~~~~~ {\rm Im} z <0 .\\
\end{array} \right.
\ee
We also have that
\[
B(t) = (1+t^2)^{-\frac{\nu}{2}} g(it).
\]
The conditions then at infinity become
\be\la{conap}
\lim_{t \rft} t^{-\nu}g(it)=1, ~~~~~~~~~~\lim_{t \ra -\infty} (-t)^{-\nu}g(it)=0.
\ee
To find the constants in (\ref{solg}) we will  satisfy the conditions at infinity (\ref{conap}) and we will
match both $g$ and $g'$ at $z=0$. That is we will ask
\be\la{mat1}
g(+i0)=g(-i0), ~~~~~~~~~~~~~g'(+i0)=g'(-i0).
\ee

We recall from \cite{AS} Section 8.1 that for $|z|>1$:
\bean
P_{\nu}(z) &  =  &  \Delta_1 z^{-\nu-1} F \left(\frac{\nu+1}{2}, \frac{\nu+2}{2}, \frac{2 \nu +3}{2}; \frac{1}{z^2} \right)
+ \Delta_2 z^{\nu} F \left(\frac{-\nu}{2}, \frac{1-\nu}{2}, \frac{1-2 \nu }{2}; \frac{1}{z^2} \right), \\
Q_{\nu}(z) &  =  &  E_1  z^{-\nu-1} F \left(\frac{\nu+2}{2}, \frac{\nu+1}{2}, \frac{2 \nu +3}{2}; \frac{1}{z^2} \right).
\eean
where,
\[
\Delta_1 = \frac{2^{-\nu-1} \pi^{-\frac{1}{2}} \Gamma(-\nu-\frac{1}{2})}{\Gamma(-\nu)},~~~~
\Delta_2 = \frac{2^{\nu} \pi^{-\frac{1}{2}} \Gamma(\nu+\frac{1}{2})}{\Gamma(1+\nu)},~~~~
E_1 =  \frac{2^{-\nu-1} \pi^{\frac{1}{2}} \Gamma(1+\nu)}{\Gamma(\frac{3}{2}+\nu)}.
\]

From the asymptotics when $t \ra \pm  \infty$,  we easily conclude
that \be\la{c11} C_1^{+} = \frac{i^{-\nu}}{\Delta_2},
~~~~~~~~~~~~~~~~~C_1^{-}=0. \ee We next see what happens near
zero. For $|z|<1$ we have that \bean P_{\nu}(z) &  =  &  B_1  F
\left(-\frac{\nu}{2}, \frac{\nu+1}{2}, \frac{1}{2}; z^2 \right)
+ B_2 z  F \left(\frac{1-\nu}{2}, \frac{2+\nu}{2}, \frac{3}{2}; z^2 \right), \\
Q_{\nu}^{\pm}(z) &  =  &  \Gamma_1    e^{\pm \frac{i \pi}{2}(-\nu-1)}
   F \left(-\frac{\nu}{2}, \frac{\nu+1}{2}, \frac{1}{2}; z^2 \right)
+ \Gamma_2 e^{\pm \frac{i \pi}{2}(-\nu)}    z   F \left(\frac{1-\nu}{2}, \frac{\nu+2}{2}, \frac{3}{2}; z^2 \right),
\eean
where the plus sign corresponds to ${\rm Im} z >0$ and the minus  to ${\rm Im} z <0$. The value of the constants
are  given by:
\bean
B_1 = \frac{\pi^{\frac{1}{2}}}{\Gamma(\frac{1-\nu}{2} )\Gamma(\frac{2+\nu}{2} )}, ~~
B_2 = \frac{-2 \pi^{\frac{1}{2}}}{\Gamma(\frac{1+\nu}{2} )\Gamma(\frac{-\nu}{2} )}, ~~
\xG_1 =\frac{\pi^{\frac{1}{2}} \Gamma(\frac{1+\nu}{2})}{2 \Gamma(1+\frac{\nu}{2} )},~~
\xG_2 =\frac{\pi^{\frac{1}{2}} \Gamma(1+\frac{\nu}{2})}{\Gamma(\frac{1+\nu}{2} )}.
\eean

 An easy calculation
shows that the matching condition (\ref{mat1}) yields
\bean
C_2^{-}  \Gamma_1  e^{ \frac{i \pi}{2}(\nu+1)} &  = &  C_1^{+} B_1  + C_2^{+} \Gamma_1  e^{ \frac{i \pi}{2}(-\nu-1)},
\\
C_2^{-}  \Gamma_2  e^{ \frac{i \pi}{2}\nu} &  = &  C_1^{+} B_2  + C_2^{+} \Gamma_2  e^{ \frac{i \pi}{2}(-\nu)},
\eean
from which it follows that
\bea\la{c22}
C_2^{+} &  = &  - \frac{ C_1^{+}}{2}  e^{ \frac{i \pi}{2}\nu} \left[ \frac{B_2}{\Gamma_2} + i \frac{B_1}{\Gamma_1} \right]
\nonumber \\
C_2^{-} &  = &   \frac{ C_1^{+}}{2}  e^{-\frac{i \pi}{2}\nu} \left[ \frac{B_2}{\Gamma_2} - i \frac{B_1}{\Gamma_1} \right].
\eea
Thus all constants in (\ref{solg}) have been computed (cf   (\ref{c11}) and (\ref{c22})), and therefore $g(z)$ is now
completely known.

The  asymptotics of $g$ for $|z| \rft$, are
\bean
g(z)   &  = &  C_1^{\pm} \Delta_2 z^{\nu} + (C_1^{\pm}\Delta_1  + C_{2}^{\pm}E_1) z^{-\nu-1} + o(|z|^{-\nu-1}),  \\
g'(z)  &  = &  C_1^{\pm} \Delta_2 \nu  z^{\nu-1} -(\nu+1) [ C_1^{\pm} \Delta_1 +C_{2}^{\pm}E_1]z^{-\nu-2} +
O(|z|^{\nu-3}),
\eean
where   the plus sign corresponds to ${\rm Im} z >0$ and the minus  to ${\rm Im} z <0$. We have that
$B(t) = (1+t^2)^{-\frac{\nu}{2}} g(it)$, whence    we  get
\[
B(t) = i^{1+\nu} C^{-}_2 E_1 (-t)^{-2 \nu -1} + o((-t)^{-2 \nu -1}),~~~~~~~~~~t \ra -\infty.
\]
 Concerning the derivative, we have for $z=it$
\[
B'(t) =-\nu t (t^2+1)^{-\frac{\nu}{2}-1}g(z) + i (1+t^2)^{-\frac{\nu}{2}}   g'(z).
\]
Whence,
\bean
B'(t) &  = &   (2 \nu +1) i^{1-\nu}(C_1^{+}\Delta_1  + C_{2}^{+}E_1)~ t^{-2 \nu -2} + o(t^{-2 \nu -2}),  ~~~~~~t \rft,  \\
B'(t) &  = &   (2 \nu +1) i^{1+\nu} C_{2}^{-}E_1 ~ (-t)^{-2 \nu -2}
 + o((-t)^{-2 \nu -2}),~~~~t \ra -\infty \ .
\eean
This completes the proof of part (ii) of the Proposition.

We next give the proof of part (i). From (\ref{tctr}) and the
asymptotics of $B(t)$ for $t \rft$, we compute \be\la{tr5} \tC =
\frac{(2 \nu+1)}{2} i^{1-2 \nu} \frac{E_1}{\Delta_2} \left( 2
\frac{\Delta_1}{E_1} -i^{\nu} \frac{B_2}{\xG_2} -i^{\nu+1}
\frac{B_1}{\xG_1} \right). \ee Using the explicit values of the
constants we calculate: \bean \frac{E_1}{\Delta_2} =  \frac{2^{-2
\nu-1}~ \pi~  \Gamma^2(1+\nu)}{ \Gamma(\frac12+\nu)
\Gamma(\frac32+\nu)}, ~~~ \frac{\Delta_1}{E_1} = \frac{ \sin(\pi
\nu)}{\pi \cos(\pi \nu)},~~~ \frac{B_2}{\xG_2} = \frac{ 2
\sin(\frac{\pi \nu}{2})}{\pi},~~~ \frac{B_1}{\xG_1} = \frac{ 2
\cos(\frac{\pi \nu}{2})}{\pi}. \eean Plugging these in (\ref{tr5})
we conclude   that  (recall that $\nu = -a/2= s-1/2$) \be\la{tr10}
\tC = \frac{2^{-2 \nu}}{\pi} \frac{\Gamma^2(1+\nu)
\Gamma(\frac12-\nu)}{\Gamma(\frac12+\nu)} = \frac{2^a
\Gamma^2(\frac{2-a}{2}) \Gamma(\frac{1+a}{2})}{  \pi
\Gamma(\frac{1-a}{2})} =\frac{2^{1-2 s}}{\pi} \frac{\Gamma^2(s
+\frac12) \Gamma(1-s)}{\Gamma(s)}. \ee

To prove part (iii) we use part (i) and we integrate the ODE (\ref{odediv2}).

By standard maximum principle arguments the solution $B(t)$ of (\ref{ode4}) subject to (\ref{bc2})
is positive and increasing.  To prove part (iv) assuming that  $a \in (-1, 0)$,  we set
  $ f(t)=(1+t^2)^{-\frac{a}{4}}B(t)$ so  that
\[
(1+t^2)f''+ 2t f'+ \frac{a(2-a)}{4} f = 0,
\]
and a similar maximum principle argument shows that $f(t)$ is also increasing. Since,
\[
f'(t)= (1+t^2)^{-\frac{a}{4}-1}  \left[(1+t^2) B' - \frac{a}{2} t B \right],
\]
we conclude that
\[
(1+t^2) B' - \frac{a}{2} t B>0, ~~~~~~~~~~t \in \R,~~~~~~~a \leq 0.
\]
Using the asymptotics of $B$, $B^{\prime}$ from part (ii) we conclude the proof of part (iv).

\finedim

Using the asymptotics of $B(t)$ from the previous  Proposition, we easily obtain  the following uniform
asymptotics for $\phi$
\begin{lemma}\la{lem33}
 Suppose  $a \in (-1,1)$  and let $\phi$ be given by
\[
\phi(x,y)   = \left\{ \begin{array}{ll}
 (y^2+d^2)^{-\frac{a}{4}} B(\frac{d}{y}),  &  ~~~~~~~  x \in \xO,~~~ y >0     \\
  (y^2+d^2)^{-\frac{a}{4}} B(-\frac{d}{y}), & ~~~~~~~ x \in \CC \xO,~~ y >0  \ ,  \\
\end{array} \right.
\]
where $B$ solves (\ref{ode4}), (\ref{bc2}).

\noindent ({\bf i})
 Then
\[
\phi(x,y)   \sim  \left\{ \begin{array}{ll}
 (y^2+d^2)^{-\frac{a}{4}},  &  ~~~~~~~  x \in \xO,~~~ y >0  \\
 y^{1-a} (y^2+d^2)^{\frac{a-2}{4}}, & ~~~~~~~ x \in \CC \xO,~~ y >0.\\
\end{array} \right.
\]
Concerning the gradient of $\phi$,  for $a \in (-1,0]$ we have
\[
|\ana \phi(x,y) |  \sim  \left\{ \begin{array}{ll}
 (y^2+d^2)^{-\frac{a+2}{4}},  &  ~~~~~~~  x \in \xO,~~~ y >0  \\
 y^{-a} (y^2+d^2)^{\frac{a-2}{4}}, & ~~~~~~~ x \in \CC \xO,~~ y >0.\\
\end{array} \right.
\]
whereas for  $a \in (0,1)$
\[
|\ana \phi(x,y) |  \sim
 y^{-a}(y^2+d^2)^{\frac{a-2}{4}},    ~~~~~~~  x \in \Ren,~~~ y >0 \ .
\]
\noindent ({\bf ii})
If     $\xO$ satisfies  $-\Delta d(x) \geq 0$ for $x \in \xO$, then  for    $a \in (-1,0)$
\[
-{\rm div}(y^a \ana \phi) \phi  \sim  \frac{y^{a}}{d (d^2+y^2)^{\frac{1+a}{2}}} (-d \Delta d) \ ,~~~~ y >0,~~~  x \in \Omega \ ,
\]
whereas  for $a=0$,
\[
-{\rm div}( \ana \phi) \phi  \sim  \frac{y}{d(d^2+y^2)} (-d \Delta d) \ ,~~~~ y >0,~~~  x \in \Omega \ .
\]
\end{lemma}

We are now ready to give the proof of Theorem \ref{th14}

\noindent
{\em Proof of Theorem \ref{th14} part (i) and (ii):} We assume that $s \in [\frac12,1)$ or
 equivalently $a \in (-1,0]$. We will use
Lemma \ref{lem31} with the test function $\phi$ given
\[
\phi(x,y)   = \left\{ \begin{array}{ll}
 (y^2+d^2)^{-\frac{a}{4}} B(\frac{d}{y}),  &  ~~~~~~~  x \in \xO,~~~ y >0     \\
  (y^2+d^2)^{-\frac{a}{4}} B(-\frac{d}{y}), & ~~~~~~~ x \in \CC \xO,~~ y >0  \ ,  \\
\end{array} \right.
\]
Using Proposition \ref{prop32} and Lemma \ref{lem33} we see that all  hypotheses of Lemma \ref{lem31} are
satisfied. In particular
 we compute
\bea\la{trco1} -\lim_{y \ra 0^+} \left( y^a \frac{\phi_{y}(x,y)
}{\phi(x,y)} \right) &  =  &
\frac{1}{d^{1-a}(x)}  \lim_{t \rft}  \left(  t^{2-a} B^{\prime}(t) \right)   \nonumber \\
&  =  & \frac{ \bar{k}_s }{d^{1-a}(x)}  \ , ~~~~~~~~~~~~~ x \in \xO \ .
\eea
We also have for $x \in \xO$ and $t=\frac{d}{y}>0$,
\bea\la{div2a}
 - {\rm div}(y^a \ana \phi) & = &-  y^a  (y^2+d^2)^{-\frac{a}{4}-1} \left[(1+t^2)^2 B''+
 (2-a)t(1+t^2) B' -\frac{a^2}{4} B  \right]   \nonumber   \\
&  & +   y^{a+1} (y^2+d^2)^{-\frac{a}{4}-1} (-\Delta d)
\left[(1+t^2) B' - \frac{a}{2} t B \right]  \nonumber   \\
& = &   y^{a+1} (y^2+d^2)^{-\frac{a}{4}-1} (-\Delta d)
\left[(1+t^2) B' - \frac{a}{2} t B \right]  \ ,
\eea
 whereas for $x \in \CC \xO$ and $t=-\frac{d}{y}<0$,  we have
\bea\la{div2b}
-{\rm div}(y^a \ana \phi) & =  & - y^a  (y^2+d^2)^{-\frac{a}{4}-1} \left[(1+t^2)^2 B''+
 (2-a)t(1+t^2) B' -\frac{a^2}{4} B  \right]   \nonumber   \\
&  & +  y^{a+1} (y^2+d^2)^{-\frac{a}{4}-1} (\Delta d)
\left[(1+t^2) B' - \frac{a}{2} t B \right] \nonumber   \\
 & =  &  y^{a+1} (y^2+d^2)^{-\frac{a}{4}-1} (\Delta d)
\left[(1+t^2) B' - \frac{a}{2} t B \right] \ .
\eea
Therefore under our assumption on $\xO$ it follows from Proposition \ref{prop32} that
\[
 - {\rm div}(y^a \ana \phi) \geq 0,  ~~~~~~~  x \in \Ren,~~~ y >0 \ .
\]
We now use  Lemma \ref{lem31} to get
\bea\la{id2b}
 \dinr y^a |\ana u|^2 \dd     \geq  &   &
\bar{k}_s  \int_{\xO}  \frac{ u^2(x,0)}{d^{1-a}(x)} dx
+ \dinr y^a |\ana u - \frac{\ana \phi}{\phi} u|^2 \dd  \nonumber \\
&  - &
 \dinr \frac{{\rm div}(y^a \ana \phi)}{\phi} u^2 \dd  \ ,
\eea
from which the trace Hardy inequality follows directly.  This relation will  also be used later on,
 in Section \ref{sec4} and  \ref{sec5} to obtain
the Sobolev term as well.

We next prove the optimality of the Hardy constant.
We will work as in section \ref{sec2}. Let
\be\la{rq2}
Q[u] := \frac{\dinr y^a |\ana u|^2 \dd}{\int_{\xO} \frac{u^2(x,0)}{d^{1-a}(x)} dx}=:\frac{N[u]}{D[u]}.
\ee
We  will show that there exists a sequence of functions $u_{\xe}$ such that
$\lim_{\xe \ra 0}Q[u_{\xe}] \leq  \tC$, and therefore $\tC$ is the best constant.

We first assume that the boundary of $\xO$  is flat in a neighborhood $U$  of a point $x_0 \in \partial \xO$.
 The  neighborhood of the point $x_0$ is assumed to  contain a ball centered
at $x_0$ with radius, say,  $3\xd$.
Locally around $x_0$ the boundary is given by $x_n=0$, whereas the interior of $\xO$ corresponds to $x_n>0$. We also
write $x=(x', x_n)$. Clearly, for $x \in \xO \cap U$  we have that $d(x) = x_n$.

We next define three suitable cutoff functions.
 Let $\psi(x') \in C_0^{\infty}(B_{\xd})$,
 where $B_{\xd} \subset \partial \xO \subset \R^{n-1}$ is the ball centered at $x_0$ with radius
$\xd$. Also the nonnegative function
 $h(x_n) \in C^{\infty}(\R)$ is such that $h(x_n) =0$  for $|x_n| \geq 2 \xd$ and $h(x_n) =1$ for $ |x_n| \leq \xd$.
We also assume that $h(x_n)$ is symmetric around $x_n=0$. Finally
 let $\chi(y) \in C_0^{\infty}(\R)$  be such that $0 \leq \chi(y) \leq 1$, and $ \chi(y) =1 $ near $y=0$.

We will use the following test function:
\be\la{ue2}
u_{\xe}(x',x_n,y)   =
 \chi(y) h(x_n) \psi(x')(y^2+x_n^2)^{-\frac{a}{4}+\frac{\xe}{4}} B(\frac{x_n}{y}), ~~~~~~~~x \in \Ren, ~~~y>0.
\ee
Using the asymptotics of $B(t)$ we easily see that
\[
u_{\xe}(x',x_n,0)   =
\left\{ \begin{array}{ll}
 h(x_n) \psi(x')x_n^{-\frac{a}{2}+\frac{\xe}{2}},    &  ~~~~~~~ x \in \xO  \\
 0, &~~~~~~~x \in \CC \xO.\\
\end{array} \right.
\]
We then compute
\be\la{de2}
D[u_{\xe}]= \int_{\R^{n-1}} \psi^2(x') dx' ~ \int_{0}^{+\infty} h^2(x_n) x_n^{-1+\xe} dx_n.
\ee
Concerning the numerator,  a straightforward calculation shows that
\bean
|\ana ( (y^2+x_n^2)^{-\frac{a}{4}+\frac{\xe}{4}} B(\frac{x_n}{y}) )|^2 = & &  \left(-\frac{a}{2}+\frac{\xe}{2} \right)^2
(y^2 + x_n^2)^{-\frac{a}{2}+\frac{\xe}{2}-1} B^2(\frac{x_n}{y}) \\
 & + &
 \frac{(x^2_n + y^2)^{1-\frac{a}{2} + \frac{\xe}{2}}}{y^4}  B^{'2}(\frac{x_n}{y}).
\eean
It is then easy to show that
\bean
N[u_{\xe}]= & &  \int_{\R^{n-1}} \psi^2(x') dx' \int_{\R} \int_{0}^{+\infty} h^2(x_n) y^a \chi^2(y) \left[
 \left(-\frac{a}{2}+\frac{\xe}{2} \right)^2
(y^2 + x_n^2)^{-\frac{a}{2}+\frac{\xe}{2}-1} B^2(\frac{x_n}{y}) \right. \\
& + & \left.
 \frac{(x_n + y^2)^{1-\frac{a}{2} + \frac{\xe}{2}}}{y^4}  B^{'2}(\frac{x_n}{y}) \right] ~  dy dx_n  +O_{\xe}(1).
\eean
To estimate the double integral above, we first break the $x_n$--integral into two pieces:  from minus infinity to zero and
 from zero to infinity. We then change variables in both pieces  by $t=x_n/y$, thus going from the  $(x_n,y)$  variables
to $(x_n,t)$. After elementary calculations we arrive at
\bean
N[u_{\xe}]= & &\int_{\R^{n-1}} \psi^2(x') dx'  ~ \int_{0}^{+\infty} h^2(x_n) x_n^{-1+\xe} dx_n \cdot \nonumber    \\
   &\cdot  & \hspace{-15mm}
  \int_{-\infty}^{+\infty}  \chi^2\left(\frac{x_n}{|t|} \right)
 \left[ \frac{(1+t^2)^{1-\frac{a}{2} + \frac{\xe}{2}}}{|t|^{\xe}} B^{'2}(t) +
 \left(-\frac{a}{2}+\frac{\xe}{2} \right)^2
 \frac{(1+t^2)^{-1-\frac{a}{2} + \frac{\xe}{2}}}{|t|^{\xe}} B^{2}(t)\right] dt  +O_{\xe}(1).
\eean
Forming the quotient we obtain
\[
Q[u_{\xe}] \leq  \int_{-\infty}^{+\infty}\left[ \frac{(1+t^2)^{1-\frac{a}{2} + \frac{\xe}{2}}}{|t|^{\xe}} B^{'2}(t) +
 \left(-\frac{a}{2}+\frac{\xe}{2} \right)^2
 \frac{(1+t^2)^{-1-\frac{a}{2} + \frac{\xe}{2}}}{|t|^{\xe}} B^{2}(t) \right] dt  +o_{\xe}(1)
\]
We finally  send $\xe$ to zero to  get
\bea\la{telosb}
\lim_{\xe \ra 0} Q[u_{\xe}] &  \leq &  \int_{-\infty}^{+\infty}\left[ (1+t^2)^{1-\frac{a}{2}} B^{'2}(t) +
 \frac{a^2}{4}
(1+t^2)^{-1-\frac{a}{2}} B^{2}(t) \right] dt  \nonumber   \\
 &  = & \tC;
\eea
the last equality follows  from Proposition \ref{prop32}(iii).

The general case where $\partial \xO$ is not  flat is treated in the same way as in section \ref{sec2}.

\finedim

\setcounter{equation}{0}
\section{Some  Weighted Hardy Inequalities}
In this section we establish some new  weighted  Hardy inequalities that will play a crucial role
in establishing trace Hardy--Sobolev--Maz'ya   inequalities.

We first prove the following:
\begin{lemma}\la{lem:abc} Let $\xO \subset \Ren$ be such that $- \Delta d(x) \geq 0$ for  $x \in \xO$.
If $A$, $B$, $\Gamma$ are constants such that $A+1>0$, $B+1>0$ and  $ 2 \Gamma< A+B+2$ then for all
 $v \in C^{\infty}_{0}(\Ren \times \R) $ there holds
\bea\la{1abc}
  \frac{(B+1)(B+A+2-2\Gamma^{+})}{B+A+2}  \dino \frac{ y^A d^{B}}{(d^2+y^2)^{\Gamma}} | v| \dd &  \leq &   \\
 \dino \frac{ y^A d^{B+1}}{(d^2+y^2)^{\Gamma}}(- \Delta d) | v| dxdy &  + &
\hspace{-2mm}
\dino \frac{ y^A d^{B+1}}{(d^2+y^2)^{\Gamma}} |\ana v| \dd \ ,   \nonumber
\eea
where $\Gamma^{+}=\max(0, \Gamma)$.
\end{lemma}

\noindent
{\em Proof:} Integrating by parts in the $x$-variables we compute
\bea\la{ipp1}
(B+1) \dino \frac{ y^A d^{B}}{(d^2+y^2)^{\Gamma}} | v| \dd =
  \dino  \frac{ y^A \ana d \cdot \ana d^{B+1}}{(d^2+y^2)^{\Gamma}} | v| dxdy \nonumber   \\
=   \dino \frac{ y^A d^{B+1}  (- \Delta d)}{(d^2+y^2)^{\Gamma}} | v| dxdy  + 2 \Gamma
 \dino \frac{ y^A d^{B+2}}{(d^2+y^2)^{\Gamma+1}} | v| \dd  \nonumber   \\
- \dino \frac{ y^A d^{B+1}}{(d^2+y^2)^{\Gamma}} \ana d \cdot \ana_x |v| \dd.
\eea
If $\Gamma \leq 0$ the result follows easily. In the sequel we consider the case $\Gamma >0$.
In the previous calculation there is no boundary term due to our assumptions.
To continue we will estimate the middle term in the right hand side above.  To this end we define
the vector field $\vec{F}$ by
\be\la{f1}
\vec{F}(x,y) := \left( \frac{ y^A d^{B+3} \ana d}{(d^2+y^2)^{\Gamma+1}},
 ~ \frac{y^{A+1} d^{B+2}}{(d^2+y^2)^{\Gamma+1}} \right).
\ee
We then have
\be\la{id1}
\dino {\rm div} \vec{F} |v| dxdy = - \dino \vec{F} \cdot \ana |v| dx dy \leq  \dino | \vec{F}| |\ana v| dxdy.
\ee
We note that because of our assumptions $A+1>0$ and $B+1>0$, there are no boundary terms in
(\ref{id1}).
Straightforward calculations show that
\be\la{div2}
{\rm div} \vec{F} = \frac{ y^A d^{B+3}  ( \Delta d)}{(d^2+y^2)^{\Gamma+1}} +
 (B+A+2-2\Gamma) \frac{ y^A d^{B+2} }{(d^2+y^2)^{\Gamma+1}},
\ee and \be\la{mod2} | \vec{F} |=  \frac{ y^A d^{B+2}
}{(d^2+y^2)^{\Gamma+1/2}} \leq \frac{ y^A d^{B+1}
}{(d^2+y^2)^{\Gamma}}. \ee From (\ref{id1})--(\ref{mod2}) we get
\bean
 (B+A+2-2\Gamma) \dino  \frac{ y^A d^{B+2} }{(d^2+y^2)^{\Gamma+1}}  |v| \dd  & &   \nonumber \\
 \leq
 \dino \frac{ y^A d^{B+3}}{(d^2+y^2)^{\Gamma+1}}(- \Delta d) | v| dxdy &  + &
\dino \frac{ y^A d^{B+1}}{(d^2+y^2)^{\Gamma}} |\ana v| \dd.
\eean
Combining the above with  (\ref{ipp1}) we conclude the proof.

\finedim

We  will also need a version of the above Lemma in case where $A+B+2=2 \xG$. In this case we have:

\begin{lemma}\la{lem:logref} Suppose that  $\xO \subset \Ren$ has finite inner radius and
 is  such that $- \Delta d(x) \geq 0$ for  $x \in \xO$.
If $A$, $B$ are constants such that $A+1>0$, $B+1>0$, then for all
 $v \in C^{\infty}_{0}(\Ren \times \R) $ there holds
\bea\la{1log}
  \frac{B+1}{A+B+3}  \dino \frac{ y^A d^{B}X^2}{(d^2+y^2)^{\frac{A+B+2}{2}}} | v| \dd &  \leq &   \\
 \dino \frac{ y^A d^{B+1}X}{(d^2+y^2)^{\frac{A+B+2}{2}}}(- \Delta d) | v| dxdy &  + &
\hspace{-2mm}
\dino \frac{ y^A d^{B+1}X}{(d^2+y^2)^{\frac{A+B+2}{2} }} |\ana v| \dd \ ,  \nonumber
\eea
where $X=X(\frac{d(x)}{R_{in}})$ and $X(t) = (1- \ln t)^{-1}$, $0<t \leq 1$.
\end{lemma}

\noindent
{\em Proof:} Integrating by parts in the $x$-variables we compute
\bea\la{4l1}
(B+1) \dino \frac{ y^A d^{B}X^2}{(d^2+y^2)^{\frac{A+B+2}{2}}} | v| \dd  + 2
\dino \frac{ y^A d^{B}X^3}{(d^2+y^2)^{\frac{A+B+2}{2}}} | v| \dd   \nonumber   \\
 \leq
 \dino \frac{ y^A d^{B+1} X^2  (- \Delta d)}{(d^2+y^2)^{\frac{A+B+2}{2}}} | v| dxdy
+ (A+B+2) \dino \frac{ y^A d^{B+2}X^2}{(d^2+y^2)^{\frac{A+B+4}{2}}} | v| \dd   \nonumber   \\
    +  \dino \frac{ y^A d^{B+1}X^2}{(d^2+y^2)^{\frac{A+B+2}{2}}} |\ana v| \dd.  ~~~~~~~~~~~~~~~
\eea

In the previous calculation there are  no boundary terms  due to our assumptions.
To continue we will estimate the middle term in the right hand side above.  To this end we define
the vector field $\vec{F}$ by
\be\la{f1log}
\vec{F}(x,y) := \left( \frac{ y^A d^{B+3} X \ana d}{(d^2+y^2)^{\frac{A+B+4}{2}} },
 ~ \frac{y^{A+1} d^{B+2} X}{(d^2+y^2)^{\frac{A+B+4}{2}}} \right).
\ee
We then have
\be\la{id1log}
\dino {\rm div} \vec{F} |v| dxdy = - \dino \vec{F} \cdot \ana |v| dx dy \leq  \dino | \vec{F}| |\ana v| dxdy.
\ee
We note that because of our assumptions $A+1>0$ and $B+1>0$, there are no boundary terms in
(\ref{id1log}).
Straightforward calculations show that
\be\la{div2log}
{\rm div} \vec{F} = \frac{ y^A d^{B+3} X  ( \Delta d)}{(d^2+y^2)^{\frac{A+B+4}{2} }} +
 \frac{ y^A d^{B+2}X^2 }{(d^2+y^2)^{\frac{A+B+4}{2}}},
\ee and \be\la{mod2log} | \vec{F} |=  \frac{ y^A d^{B+2} X
}{(d^2+y^2)^{\frac{A+B+3}{2} }} \leq \frac{ y^A d^{B+1}X
}{(d^2+y^2)^{\frac{A+B+2}{2}}}. \ee From
(\ref{id1log})--(\ref{mod2log}) we get \bean
  \dino  \frac{ y^A d^{B+2}X^2 }{(d^2+y^2)^{\frac{A+B+4}{2} }}  |v| \dd  & &   \nonumber \\
 \leq
 \dino \frac{ y^A d^{B+3}X }{(d^2+y^2)^{\frac{A+B+4}{2} }}(- \Delta d) | v| dxdy &  + &
\dino \frac{ y^A d^{B+1}X }{(d^2+y^2)^{\frac{A+B+2}{2}}} |\ana v| \dd.
\eean
Combining the above with  (\ref{4l1}) we conclude the proof.

\finedim

Without imposing any geometric assumption on $\xO$ we have the following result that will also be used later on.

\begin{lemma}\la{lem:gen} Let $\xO \subset \Ren$.
If $A$, $B$, $\Gamma$ are constants such that $A+1>0$, $B+1>0$ and  $ 2 \Gamma< A+B+2$, then there exist
positive constants $c_1$ and $c_2$ such that  for all
 $v \in C^{\infty}_{0}(\Ren \times \R) $ there holds
\bea\la{4.gen}
   \dino \frac{ y^A d^{B}}{(d^2+y^2)^{\Gamma}} | v| \dd
 ~~~~~~~~~~~~~~~~~~~~~~~~~~~~~~~~~~~~~~~~~~~~~~~~~~ ~~~~~~~~~   ~~ ~~~~~~~~~ ~~~ ~~~~~~   ~~ ~~     \\
 ~~~~~~~~~~  \leq  c_1  \dino \frac{ y^{A} d^{B+1}}{(d^2+y^2)^{\Gamma}} |\nabla v| dxdy   +
  c_2
\dino \frac{ y^{A} d^{B+1}}{(d^2+y^2)^{\Gamma}} | v| \dd \ .   \nonumber
\eea

\end{lemma}

\noindent
{\em Proof:} Here we will use the fact that $\partial \xO$ is uniformly Lipschitz. Let  $\{U_i\}$ be a  covering of
$\xO_{\xe} = \{x \in \xO: {\rm dist} (x, \partial \xO)<\xe \}$ and let $\phi_i$ be a partition of unity subordinate
to the covering $\{U_i\}$. We then have
\[
 \int_{0}^{+\infty}  \int_{\xO_{\xe}} \frac {y^{A} d^{B}}{(d^2+y^2)^{\Gamma}} | v| \dd
\leq \sum_{i=1}^{+\infty}  \int_{0}^{+\infty}  \int_{\xO_{\xe}}
\frac {y^{A} d^{B}}{(d^2+y^2)^{\Gamma}} |\phi_i v| dd.
\]
In each $U_i$ we straighten the boundary and use the equivalence of the distance function to the regularized distance
as well as to the difference $x_n - f_i(x')$ (see \cite{St} section 3.2, or \cite{L} section 12.2) and obtain
\[
\int_{0}^{+\infty}  \int_{\xO_{\xe}} \frac {y^{A}
d^{B}}{(d^2+y^2)^{\Gamma}} |\phi_i v| \dd \leq C
\int_{0}^{+\infty} \int_{\R^n_{+}} \frac {y^{A} t^{B}}
{(t^2+y^2)^{\Gamma}} |\tilde{\phi}_i \tilde{v}| \dd  \ ,
\]
for some constant $C$ independent of $i$. We next use Lemma
\ref{lem:abc} to estimate the right hand side of this, thus
obtaining \bean \int_{0}^{+\infty}  \int_{\R^n_{+}} \frac {y^{A}
t^{B}} {(t^2+y^2)^{\Gamma}} |\tilde{\phi}_i \tilde{v}| \dd \leq C
\int_{0}^{+\infty}  \int_{\R^n_{+}} \frac {y^{A} t^{B+1}} {(t^2+y^2)^{\Gamma}} |\nabla(  \tilde{\phi}_i \tilde{v})| \dd \\
\leq C \int_{0}^{+\infty}  \int_{\R^n_{+}} \frac {y^{A} t^{B+1}}
{(t^2+y^2)^{\Gamma}}  \tilde{\phi}_i|  \nabla  \tilde{v}| \dd + C
\int_{0}^{+\infty}  \int_{\R^n_{+}} \frac {y^{A} t^{B+1}}
{(t^2+y^2)^{\Gamma}}  | \nabla \tilde{\phi}_i||\tilde{v}| \dd
\eean Hence, returning to our original variables we have that
\bean
\int_{0}^{+\infty}  \int_{\xO_{\xe}} \frac {y^{A} d^{B}}{(d^2+y^2)^{\Gamma}} |\phi_i v| \dd \hspace{6cm}  \\
 \leq
C \int_{0}^{+\infty}  \int_{\xO_{\xe}} \frac {y^{A} d^{B+1}}
{(d^2+y^2)^{\Gamma}} \phi_i|  \nabla  v| \dd + C
\int_{0}^{+\infty}  \int_{\xO_{\xe}} \frac {y^{A} d^{B+1}}
{(d^2+y^2)^{\Gamma}}  | \nabla \phi_i||v| \dd  \ . \eean Summing
over $i$ we get that \bean
 \int_{0}^{+\infty}  \int_{\xO_{\xe}} \frac {y^{A} d^{B}}{(d^2+y^2)^{\Gamma}} | v| \dd  \hspace{6cm}  \\  \leq
C_1 \int_{0}^{+\infty}  \int_{\xO_{\xe}} \frac {y^{A}
d^{B+1}}{(d^2+y^2)^{\Gamma}} |\nabla v| \dd + C_2
\int_{0}^{+\infty}  \int_{\xO_{\xe}} \frac {y^{A}
d^{B+1}}{(d^2+y^2)^{\Gamma}} | v| \dd   \ . \eean The result then
follows easily.

\finedim

When working in the complement of $\xO$ we have the following surprising result:

\begin{lemma}\la{lem:abc22} Let $\xO \subset \Ren$.
If $A$, $B$, $\Gamma$ are constants such that $A+1>0$, $B+1>0$ and  $ 2 \Gamma< A+B+2$ then for all
 $v \in C^{\infty}_{0}(\Ren \times \R) $ there holds
\bea\la{1abc22}
  (A+1)(A+B+2-2\Gamma^{+})  \dinoC \frac{ y^A d^{B}}{(d^2+y^2)^{\Gamma}} | v| \dd   \leq
 ~~~~~~~~~~~~~~~~~~~~~~~~~~~~~~~~~~~~~~~~~~~~~~~~~~       \\
 2 \xG^{+}  \dinoC \frac{ y^{A+2} d^{B+1}}{(d^2+y^2)^{\Gamma+1}}(- \Delta d) | v| dxdy   +
  (A+B+2)
\dinoC \frac{ y^{A+1} d^{B}}{(d^2+y^2)^{\Gamma}} |\ana v| \dd \ ,   \nonumber
\eea
where $\Gamma^{+}=\max(0, \Gamma)$.
\end{lemma}
We note that no assumption on the sign of $- \Delta d$ is required.

\noindent
{\em Proof:} Integrating  by parts in the $y$-variable we compute
\bea\la{ipp22}
(A+1) \dinoC \frac{ y^A d^{B}}{(d^2+y^2)^{\Gamma}} | v| \dd  \leq  2 \xG
 \dinoC \frac{ y^{A+2} d^{B}}{(d^2+y^2)^{\Gamma+1}} | v| \dd   \nonumber    \\
 + \dinoC \frac{ y^{A+1} d^{B}}{(d^2+y^2)^{\Gamma}} | \ana v| \dd.
\eea
If $\Gamma \leq 0$ the result follows easily. In the sequel we consider the case $\Gamma >0$.
In the previous calculation there is no boundary term due to our assumptions.
To continue we will estimate the first term in the right hand side above.  To this end we define
the vector field $\vec{F}$ by
\be\la{f122}
\vec{F}(x,y) := \left( \frac{ y^{A+2} d^{B+3} \ana d}{(d^2+y^2)^{\Gamma+1}},
 ~ \frac{y^{A+3} d^{B}}{(d^2+y^2)^{\Gamma+1}} \right).
\ee
We then have
\be\la{id122}
\dinoC {\rm div} \vec{F} |v| dxdy = - \dinoC \vec{F} \cdot \ana |v| dx dy \leq  \dinoC | \vec{F}| |\ana v| dxdy.
\ee
We note that because of our assumptions $A+1>0$ and $B+1>0$, there are no boundary terms in
(\ref{id122}).
Straightforward calculations show that
\be\la{div222}
{\rm div} \vec{F} = \frac{ y^{A+2} d^{B+1}  ( \Delta d)}{(d^2+y^2)^{\Gamma+1}} +
 (A+B+2-2\Gamma) \frac{ y^{A+2} d^{B} }{(d^2+y^2)^{\Gamma+1}},
\ee
and
\be\la{mod222}
| \vec{F} |=  \frac{ y^{A+2} d^{B} }{(d^2+y^2)^{\Gamma+1/2}} \leq \frac{ y^{A+1} d^{B} }{(d^2+y^2)^{\Gamma}}.
\ee
Combining the above we conclude the proof.
Again,we note that in all integrations by  parts  there are  no boundary terms due to our assumptions.

\finedim

As a consequence of Lemma \ref{lem:abc} we have:
\begin{lemma}\la{lem:l2}
 Let $\xO \subset \Ren$ be such that $- \Delta d(x) \geq 0$, for  $x \in \xO$ and
 $w \in C^{1}_{0}(\Ren \times \R) $.
If $A$, $B$, $\Gamma$ are constants such that $A+1>0$, $B+1>0$, and  $ 2 \Gamma< A+B+2$,
then,
\bea\la{1abc2}
  \frac{(B+1)^2(B+A+2-2\Gamma^{+})^2}{4(B+A+2)^2}  \dino \frac{ y^A d^{B}}{(d^2+y^2)^{\Gamma}}  w^2 \dd &  \leq &   \\
 \frac{(B+1)(B+A+2-2\Gamma^{+})}{2(B+A+2)}  \dino \frac{ y^A d^{B+1}}{(d^2+y^2)^{\Gamma}}(- \Delta d) w^2 dxdy &  + &
\hspace{-2mm}
\dino \frac{ y^A d^{B+2}}{(d^2+y^2)^{\Gamma}} |\ana w|^2 \dd \ ,  \nonumber
\eea
where $\Gamma^{+}=\max(0, \Gamma)$.
\end{lemma}
\noindent
{\em Proof:} We apply Lemma \ref{lem:abc} to $v=w^2$.
To conclude we use Young's inequality in the last term of the right hand side. We omit the details.

\finedim

In the case where $A+B+2=2\xG$ the $L^2$ analogue of Lemma \ref{lem:logref} reads:

\begin{lemma}\la{lem:l2log}
  Suppose that  $\xO \subset \Ren$ has finite inner radius and
 is  such that $- \Delta d(x) \geq 0$ for  $x \in \xO$.
If $A$, $B$ are constants such that $A+1>0$, $B+1>0$, then for all
 $w \in C^{\infty}_{0}(\Ren \times \R) $ there holds
\bea\la{1logl2}
 \left( \frac{B+1}{2(A+B+3)} \right)^2
  \dino \frac{ y^A d^{B}X^2}{(d^2+y^2)^{\frac{A+B+2}{2}}}  w^2 \dd &  \leq &   \\
 \frac{B+1}{2(A+B+3)}   \dino \frac{ y^A d^{B+1}X}{(d^2+y^2)^{\frac{A+B+2}{2}}}(- \Delta d)  w^2 dxdy &  + &
\hspace{-2mm}
\dino \frac{ y^A d^{B+2}}{(d^2+y^2)^{\frac{A+B+2}{2} }} |\ana w|^2 \dd \ ,  \nonumber
\eea
where $X=X(\frac{d(x)}{R_{in}})$ and $X(t) = (1- \ln t)^{-1}$, $0<t \leq 1$.

\end{lemma}
\noindent
{\em Proof:} We apply Lemma \ref{lem:logref} to $v=w^2$.
To conclude we use Young's inequality in the last term of the right hand side. We omit the details.

\finedim

In the case of half space a more delicate result is needed. More precisely we have:

\begin{lemma}\la{lem:r2}
Let $v \in C_{0}^{\infty}(\R^n \times \R)$. If  $0<A \leq  \frac12$, $B+1>0$, and  $ 2 \Gamma< A+B+2$,  then
the following inequality holds true:
\bea\la{4.6}
 c_0  \int_{0}^{
 +\infty} \int_{\R^n_{+}}
\frac{y^{-A}x_n^{B}}{(x_n^2+y^2)^{\Gamma -A}} |v| dx dy  \leq
\int_{0}^{+\infty}
\int_{\R^n_{+}}\frac{y^{A}x_n^{1+B}}{(x_n^2+y^2)^{\Gamma}}|\ana v|
dx dy \ , \eea where
\[
c_0=  \frac{A(B+1)(B+A+2-2\Gamma^{+})}{(A+B+2)(A+2B+2)-2\Gamma^{+} (B+1)} \ .
\]
The same result holds true if we replace $\R^n_{+}$ by $\R^n_{-}$ with $|x_n|$ in the place of $x_n$.
\end{lemma}

\noindent {\em Proof:} We will use polar coordinates, $x_n=r \cos \theta$, $y=r \sin \theta$.
 We first establish the following inequality for the
angular derivative.
\bea\la{4.10}
A \int_{0}^{\frac{\pi}{2}} ( \sin \theta)^{-A} ( \cos \theta)^{B}|v|   d \theta
&\leq &     (1+A+B) \int_{0}^{\frac{\pi}{2}}( \sin \theta)^{1+A}  ( \cos \theta)^{B}|v| d  \theta \nonumber   \\
& +&  \int_{0}^{\frac{\pi}{2}}( \sin \theta)^{A}  ( \cos \theta)^{1+B}|v_{\theta}|   d \theta \ .
\eea
We have
\bean
\frac{d}{d \theta} (( \sin \theta)^{A}  ( \cos \theta)^{1+B}) &  = &  A ( \sin \theta)^{A-1}  ( \cos \theta)^{2+B} -
(1+B) ( \sin \theta)^{A+1}  ( \cos \theta)^{B}   \\
 &  = &  A ( \sin \theta)^{A-1}  ( \cos \theta)^{B} -
(1+A+B) ( \sin \theta)^{A+1}  ( \cos \theta)^{B} \ ,
\eean
therefore an integration by parts gives:
\bean
A \int_{0}^{\frac{\pi}{2}} ( \sin \theta)^{A-1} ( \cos \theta)^{B}|v|   d \theta & \leq &
(1+A+B)
 \int_{0}^{\frac{\pi}{2}}( \sin \theta)^{1+A}  ( \cos \theta)^{B}|v| d  \theta \nonumber   \\
& +&  \int_{0}^{\frac{\pi}{2}}( \sin \theta)^{A}  ( \cos \theta)^{1+B}|v_{\theta}|   d \theta \ .
\eean
Since $A \leq \frac{1}{2}$ we also have that $ ( \sin \theta)^{-A} \leq ( \sin \theta)^{A-1}$ and (\ref{4.10})
follows.

We next multiply (\ref{4.10}) by $r^{A+B+1-2\xG}$ and then
integrate over $(0, \infty)$ to conclude: \bea A
\int_{0}^{+\infty} \int_{0}^{+\infty}
\frac{y^{-A}x_n^{B}}{(x_n^2+y^2)^{\Gamma -A}} |v| dx_n dy &   \leq
& (1+A+B)   \int_{0}^{+\infty} \int_{0}^{+\infty}
\frac{y^{1+A}x_n^{B}}{(x_n^2+y^2)^{\Gamma +\frac12}} |v| dx_n dy  \nonumber    \\
&  + &
   \int_{0}^{+\infty} \int_{0}^{+\infty}\frac{y^{A}x_n^{1+B}}{(x_n^2+y^2)^{\Gamma}}|\ana v| dx_n dy  \nonumber   \\
 &   \leq  &
(1+A+B)   \int_{0}^{+\infty} \int_{0}^{+\infty}
\frac{y^{A}x_n^{B}}{(x_n^2+y^2)^{\Gamma}} |v| dx_n dy  \nonumber    \\
&  + &
   \int_{0}^{+\infty} \int_{0}^{+\infty}\frac{y^{A}x_n^{1+B}}{(x_n^2+y^2)^{\Gamma}}|\ana v| dx_n dy  \ .
\eea
We next estimate the first term in the right hand side by using   Lemma \ref{lem:abc}, that is,
\[
 \frac{(B+1)(B+A+2-2\Gamma^{+})}{B+A+2}  \int_{0}^{+\infty} \int_{0}^{+\infty}
  \frac{ y^A x_n^{B}}{(x_n^2+y^2)^{\Gamma}} | v| dx_n dy   \leq
 \int_{0}^{+\infty} \int_{0}^{+\infty}     \frac{ y^A x_n^{B+1}}{(x_n^2+y^2)^{\Gamma}} |\ana v| dx_n dy \ .
\]
A further integration in the other variables completes the proof.

\finedim

\setcounter{equation}{0}
\section{Half  Space, Trace Hardy   \&  Trace Hardy--Sobolev--Maz'ya   Inequalities}\la{sec4}
Here we will prove the trace Hardy and trace  Hardy--Sobolev--Maz'ya inequalities appearing in Theorems
 \ref{th12} and \ref{th16}.
We start with  the trace Hardy inequalities.

\subsection{Half  Space,  Trace Hardy  I   \& II}
In this subsection we will provide the proof of  the trace Hardy inequalities appearing in Theorems
\ref{th12} and \ref{th16}.

{\em Proof of Theorem  \ref{th12} part (i) and (ii):} The case where $s \in [\frac12,1)$
is contained in Theorem \ref{th11}. We next consider the case $s  \in (0,\frac12)$ or equivalently $a \in (0,1)$.

We will use the notation $x=(x^{\prime}, x_n) \in \Ren_{+}$ with $x_n>0$.  We will use Lemma  \ref{lem21}
with   the test function $\phi$ given by
\[
\phi(x,y)= x_n^{-\frac{a}{2}} A \left( \frac{y}{x_n} \right), ~~~~~~~~
 y >0,~~~ x_n>0, x \in \Ren_{+}  \ ,
\]
where $A$ solves (\ref{ode1}), (\ref{2bc}).  Using Proposition \ref{prop22}   and Lemma  \ref{lemasf}  we see
that all hypotheses of Lemma \ref{lem21} are satisfied. In particular,
 for $t = \frac{y}{x_n}$ we compute, for $x \in \Ren_{+}$,
\[
-\lim_{y \ra 0^+} \left( y^a \frac{\phi_{y}(x,y) }{\phi(x,y)} \right)   =
   \frac{ \bar{d}_s }{x_n^{1-a}}  \ .
 \]
We also have
\[
- {\rm div}(y^a \ana \phi)  = 0, ~~~~~  y >0,~~~~~ x \in \Ren_{+}  \ .
\]
From Lemma \ref{lem21} we get \be\la{id2hs}
 \int_0^{+\infty} \int_{\R^n_{+}}  y^a |\ana u|^2 \dd     \geq
\bar{d}_s  \int_{\R^n_{+}}  \frac{ u^2(x,0)}{x_n^{1-a}} dx +
\int_0^{+\infty} \int_{\R^n_{+}}   y^a |\ana u - \frac{\ana
\phi}{\phi} u|^2 \dd \ee from which the trace Hardy inequality
follows directly.  This relation will be used later on,
% in Sections \ref{sec4}  and \ref{sec5}
 to obtain the Sobolev term as well.

The optimality of $ \bar{d}_s$ follows by the same test functions given by (\ref{ue})
  as in the flat case of Theorem \ref{th11}.
The fact that $a$  covers the full interval $(-1,1)$ does not affect the calculations leading to
(\ref{telos}).

\finedim

{\em Proof of Theorem \ref{th16} part (i):}
 The case where $s \in [\frac12,1)$
is contained in Theorem \ref{th14}. We next consider the case $s  \in (0,\frac12)$ or equivalently $a \in (0,1)$.
 We will use
Lemma \ref{lem31} with the test function $\phi$ given
\[
\phi(x,y)   =
 (y^2+x_n^2)^{-\frac{a}{4}} B(\frac{x_n}{y}),    ~~~~~~ y >0, ~~x_n \in \R \ .
\]
Using Proposition \ref{prop32} and Lemma \ref{lem33} we see that all  hypotheses of Lemma \ref{lem31} are
satisfied. In particular
 we compute
\[
-\lim_{y \ra 0^+} \left( y^a \frac{\phi_{y}(x,y) }{\phi(x,y)} \right)
  =   \frac{ \bar{k}_s }{x_n^{1-a}}  \ , ~~~~~~~~~~~~~ x_n>0  \ .
\]
An easy calculation shows that
\[
 - {\rm div}(y^a \ana \phi)  = 0, ~~~~~~~~x \in \R^n,~~~y>0 \ .
\]

We now use  Lemma \ref{lem31} to get
\bea\la{id2bsob}
 \dinr y^a |\ana u|^2 \dd     \geq
\bar{k}_s  \int_{\R^n_{+}}  \frac{ u^2(x,0)}{x_n^{1-a}} dx
+ \dinr y^a |\ana u - \frac{\ana \phi}{\phi} u|^2 \dd
\eea
from which the trace Hardy inequality follows directly.  This relation will  also be used later on,
 to obtain  the Sobolev term as well.

The optimality of $ \bar{k}_s$ follows by the same test functions given by (\ref{ue2})
  as in the flat case of Theorem \ref{th14}.
The fact that $a$  covers the full interval $(-1,1)$ does not affect the calculations leading to
(\ref{telosb}).

\finedim

\subsection{Half  Space,  Trace Hardy--Sobolev--Maz'ya I   \& II }

Here we will give the proof of the trace  Hardy--Sobolev--Maz'ya inequalities of  Theorems
\ref{th12} and \ref{th16}. We will first establish  different trace Hardy--Sobolev--Maz'ya inequalities
where only the Hardy term appears in the trace, and which are of independent interest.

\begin{theorem}\la{th45}
Let $0<s<1$ and  $n \geq 2$.
 There exists a positive constant $c$ such that
 for all   $u \in C^{\infty}_{0}(\R^n_{+} \times \R) $ there holds
\be\la{4.13}
 \int_{0}^{+\infty} \int_{\R^n_{+}}  y^{1-2s} |\ana_{(x,y)} u(x,y)|^2 \dd  \geq  \bar{d}_{s}
 \int_{\R^n_{+}} \frac{u^2(x,0)}{x_n^{2s}}dx + c
 \left(\int_{0}^{+\infty}  \int_{\R^n_{+}} | u(x,y)|^{\frac{2(n+1)}{n-2s}} dx dy \right)^{\frac{n-2s}{n+1}}  \ .
\ee
with
\be\la{4.13d}
 \bar{d}_{s} :=  \frac{2 \Gamma \left(1-s \right) \Gamma^2 \left(\frac{3+2s}{4} \right)}
{\Gamma^2 \left(\frac{3-2s}{4} \right) \Gamma \left( s \right)} \ .
\ee
\end{theorem}

\noindent {\em Proof of Theorem \ref{th45}:} From the proof of Theorem \ref{th12} we recall the inequality
(\ref{id2hs}), that is
\be\la{4hs}
 \int_0^{+\infty} \int_{\R^n_{+}}  y^a |\ana u|^2 \dd     \geq
\bar{d}_s  \int_{\R^n_{+}}  \frac{ u^2(x,0)}{x_n^{1-a}} dx +
\int_0^{+\infty} \int_{\R^n_{+}}   y^a |\ana u - \frac{\ana
\phi}{\phi} u|^2 \dd  \ , \ee where
 $\phi$ is  given by
\[
\phi(x,y)= x_n^{-\frac{a}{2}} A \left( \frac{y}{x_n} \right), ~~~~~~~~
 y >0,~~~~ x_n>0,~~~ x \in \Ren_{+}  \ ,
\]
and  $A$ solves (\ref{ode1}), (\ref{2bc}).

The result will follow after  establishing  the following inequality:
\be\la{4.28}
 \int_0^{+\infty} \int_{\R^n_{+}}   y^a |\ana u - \frac{\ana \phi}{\phi} u|^2 \dd
\geq c \left(\int_0^{+\infty} \int_{\R^n_{+}}
|u|^{\frac{2(n+1)}{n+a-1}} \dd \right)^{\frac{n+a-1}{(n+1)}} \ .
\ee

To this end we start with  the inequality,  see \cite{Maz}, Theorem 1, section 2.1.6,
\[
\int_0^{+\infty} \int_{\R^n_{+}}  y^{\frac{a}{2}} |\ana u| \dd
\geq c \left(\int_0^{+\infty} \int_{\R^n_{+}}
|u(x,y)|^{\frac{2(n+1)}{2n+a}} \dd \right)^{\frac{2n+a}{2(n+1)}} \
,~~~~~~~ u \in C^{\infty}_{0}(\R^n_{+} \times \R) \ ,
\]
with the choice $u = \phi^{\frac{2n+a}{n+a-1}} v$. Hence we obtain
\bea\la{4.30} \int_0^{+\infty} \int_{\R^n_{+}}   y^{\frac{a}{2}}
\phi^{\frac{2n+a}{n+a-1}}  |\ana v| \dd + \frac{2n+a}{n+a-1}
\int_0^{+\infty} \int_{\R^n_{+}}  y^{\frac{a}{2}}
\phi^{\frac{n+1}{n+a-1}} |\ana \phi|
|v| \dd     \nonumber         \\
 \geq
 c \left( \int_0^{+\infty} \int_{\R^n_{+}} |\phi^{\frac{2n+a}{n+a-1}} v|^{\frac{2(n+1)}{2n+a}}
\dd \right)^{\frac{2n+a}{2(n+1)}}   \ .
\eea
Next we will control the second term of the LHS by the first term of the LHS. To this end we consider
two cases. Suppose first that $s \in [\frac12,1)$ that is $a \in (-1,0]$. Using the asymptotics of
Lemma  \ref{lemasf} we get that
\[
 y^{\frac{a}{2}}  \phi^{\frac{n+1}{n+a-1}} |\ana \phi|
\sim
 \frac{ y^{\frac{a}{2}}  x_n^{\frac{n+1}{n+a-1}}}{(x_n^2+y^2)^{\frac{(2+a)(2n+a)}{4(n+a-1)}}} \ ,
\]
whereas,
\be\la{4.34}
 y^{\frac{a}{2}}  \phi^{\frac{2n+a}{n+a-1}} \sim
 \frac{ y^{\frac{a}{2}}  x_n^{\frac{2n+a}{n+a-1}}}{(x_n^2+y^2)^{\frac{(2+a)(2n+a)}{4(n+a-1)}}} \ .
\ee
The sought for estimate then is a consequence of Lemma \ref{lem:abc} with the choice:
$A= \frac{a}{2}$, $B = \frac{n+1}{n+a-1}$ and $\xG = \frac{(2+a)(2n+a)}{4(n+a-1)}$ taking into
account that
\[
A+B+2-2\xG = \frac{(2-a)(n-1)}{2(n+a-1)} >0  \ .
\]
We next consider the case $a \in (0,1)$. Using again the asymptotics of Lemma \ref{lemasf} this time we have that
\[
 y^{\frac{a}{2}}  \phi^{\frac{n+1}{n+a-1}} |\ana \phi|
\sim
 \frac{ y^{-\frac{a}{2}}  x_n^{\frac{n+1}{n+a-1}}}{(x_n^2+y^2)^{\frac{(2+a)(n+1)}{4(n+a-1)}+\frac{2-a}{4}}} \ ,
\]
whereas, (\ref{4.34}) remains the same.
The sought for estimate now  is a consequence of Lemma \ref{lem:r2} with the choice
$A= \frac{a}{2}$, $B = \frac{n+1}{n+a-1}$ and $\xG = \frac{(2+a)(2n+a)}{4(n+a-1)}$ taking into
account that
\[
A+B+2-2\xG = \frac{(2-a)(n-1)}{2(n+a-1)} >0  \ .
\]

 Therefore  for any $a \in (-1,1)$ we arrive at:
\bea\la{4.36} \int_0^{+\infty} \int_{\R^n_{+}}   y^{\frac{a}{2}}
\phi^{\frac{2n+a}{n+a-1}}  |\ana v|
 \geq
 c \left( \int_0^{+\infty} \int_{\R^n_{+}} |\phi^{\frac{2n+a}{n+a-1}} v|^{\frac{2(n+1)}{2n+a}}
\dd \right)^{\frac{2n+a}{2(n+1)}}   \ . \eea To continue we next
set in (\ref{4.36})  $v = |w|^{\frac{2n+a}{n+a-1}}$ and apply
Schwartz inequality in the LHS to conclude after a simplification
\be\la{4.38} \int_0^{+\infty} \int_{\R^n_{+}}  y^{a} \phi^2 |\ana
w|^2 \dd \geq c \left(\int_0^{+\infty} \int_{\R^n_{+}} |\phi
w|^{\frac{2(n+1)}{n+a-1}} \dd \right)^{\frac{n+a-1}{n+1}} \ , \ee
which is equivalent to (\ref{4.28}).

\finedim

{\em Proof of Theorem \ref{th12} part (iii):} Our starting point now is the following weighted  trace Sobolev
inequality,   see \cite{Maz}, Theorem 1, section 2.1.6,
\[
\int_0^{+\infty} \int_{\R^n_{+}}  y^{\frac{a}{2}} |\ana u| \dd
\geq c \left( \int_{\R^n_{+}} |u(x,0)|^{\frac{2n}{2n+a}} dx
\right)^{\frac{2n+a}{2n}} \ ,~~~~~~~ u \in C^{\infty}_{0}(\R^n_{+}
\times \R) \ .
\]
Again we set $u = \phi^{\frac{2n+a}{n+a-1}} v$, to obtain the
analogue of (\ref{4.30}). \bea\la{4.40} \int_0^{+\infty}
\int_{\R^n_{+}}   y^{\frac{a}{2}}  \phi^{\frac{2n+a}{n+a-1}} |\ana
v| \dd + \frac{2n+a}{n+a-1} \int_0^{+\infty} \int_{\R^n_{+}}
y^{\frac{a}{2}}  \phi^{\frac{n+1}{n+a-1}} |\ana \phi|
|v| \dd     \nonumber         \\
 \geq
 c \left( \int_{\R^n_{+}} |\phi^{\frac{2n+a}{n+a-1}}(x,0) v(x,0)|^{\frac{2n}{2n+a}}
dx \right)^{\frac{2n+a}{2n}}   \ .
\eea
As in the proof of Theorem \ref{th45}
 we  control the second term of the LHS by the first term of the LHS to arrive at
\[
\int_0^{+\infty} \int_{\R^n_{+}}   y^{\frac{a}{2}}
\phi^{\frac{2n+a}{n+a-1}}  |\ana v| \dd \geq c \left(
\int_{\R^n_{+}} |\phi^{\frac{2n+a}{n+a-1}}(x,0)
v(x,0)|^{\frac{2n}{2n+a}} dx \right)^{\frac{2n+a}{2n}}   \ .
\]
Again, we set  $v = |w|^{\frac{2n+a}{n+a-1}}$ and apply Schwartz
inequality in the LHS to arrive at \bean \left( \int_0^{+\infty}
\int_{\R^n_{+}}  y^{a} \phi^2 |\ana w|^2 \dd \right)^{\frac12}
\left(\int_0^{+\infty} \int_{\R^n_{+}} |\phi
w|^{\frac{2(n+1)}{n+a-1}} \dd \right)^{\frac12} \geq  c \left(
\int_{\R^n_{+}} |(\phi w)(x,0) |^{\frac{2n}{n+a-1}} dx
\right)^{\frac{2n+a}{2n}}   \ . \eean We next use (\ref{4.38}) to
conclude after a simplification
\[
 \int_0^{+\infty} \int_{\R^n_{+}}  y^{a} \phi^2 |\ana w|^2 \dd
\geq c
\left( \int_{\R^n_{+}} |(\phi w)(x,0) |^{\frac{2n}{n+a-1}}
dx \right)^{\frac{n+a-1}{n}}   \ ,
\]
which is equivalent to
\[
 \int_0^{+\infty} \int_{\R^n_{+}}   y^a |\ana u - \frac{\ana \phi}{\phi} u|^2 \dd
\geq c
\left( \int_{\R^n_{+}} |u(x,0) |^{\frac{2n}{n+a-1}}
dx \right)^{\frac{n+a-1}{n}}   \ .
\]
Combining this with inequality
(\ref{id2hs}) we conclude the proof.

\finedim

We next present a preliminary result which will play an important role towards establishing the
Hardy--Sobolev--Maz'ya II of Theorem \ref{th16}.

\begin{theorem}\la{th46}
Let $0<s<1$ and  $n \geq 2$.
   There exists a positive constant $c$, such that
for all   $u \in C^{\infty}_{0}(\Ren \times \R) $
with  $u(x,0)=0$, $x \in \Ren_{-}$, there holds
\be\la{4.C1s}
 \dinr y^{1-2s} |\ana_{(x,y)} u(x,y)|^2 \dd  \geq  \bar{k}_{s}
 \int_{\Ren_{+}} \frac{u^2(x,0)}{x_n^{2s}}dx + c
\left(\int_{0}^{+\infty}  \int_{\R^n} |
u(x,y)|^{\frac{2(n+1)}{n-2s}} dx dy \right)^{\frac{n-2s}{n+1}}
  \ ,
\ee
where
\[
  \bar{k}_{s} :=  \frac{ 2^{1-2 s}  \Gamma^2(s +\frac12) \Gamma(1-s)}{\pi \Gamma(s)}  \ ,
\]
is the best constant in (\ref{4.C1s}).
\end{theorem}

\noindent {\em Proof:}
 From the proof of Theorem \ref{th14} we recall the inequality    (\ref{id2b}),
 that is
\be\la{44hs}
 \int_0^{+\infty} \int_{\R^n}  y^a |\ana u|^2 \dd     \geq
\bar{k}_s  \int_{\R^n_{+}}  \frac{ u^2(x,0)}{x_n^{1-a}} dx +
\int_0^{+\infty} \int_{\R^n}   y^a |\ana u - \frac{\ana
\phi}{\phi} u|^2 \dd  \ , \ee where
 $\phi$ is  given by
\[
\phi(x,y)   =
 (y^2+x_n^2)^{-\frac{a}{4}} B(\frac{x_n}{y}),    ~~~~~~ y >0, ~~x_n \in \R \ ,
\]
and $B$ solves (\ref{ode4}), (\ref{bc2}).

Again, the result will follow after  establishing  the following inequality:
\be\la{4.48}
 \int_0^{+\infty} \int_{\R^n}   y^a |\ana u - \frac{\ana \phi}{\phi} u|^2 \dd
   \geq c
\left(\int_0^{+\infty} \int_{\R^n} |u|^{\frac{2(n+1)}{n+a-1}} \dd
\right)^{\frac{n+a-1}{n+1}} \ . \ee To this end we start with
the inequality,  see \cite{Maz}, Theorem 1, section 2.1.6,
\[
\int_0^{+\infty} \int_{\R^n}  y^{\frac{a}{2}} |\ana u| \dd \geq c
\left(\int_0^{+\infty} \int_{\R^n} |u(x,y)|^{\frac{2(n+1)}{2n+a}}
\dd \right)^{\frac{2n+a}{2(n+1)}} \ ,~~~~~~~ u \in
C^{\infty}_{0}(\R^n \times \R) \ ,
\]
with the choice $u = \phi^{\frac{2n+a}{n+a-1}} v$. Hence we obtain
\bea\la{4.60} \int_0^{+\infty} \int_{\R^n}   y^{\frac{a}{2}}
\phi^{\frac{2n+a}{n+a-1}}  |\ana v| \dd + \frac{2n+a}{n+a-1}
\int_0^{+\infty} \int_{\R^n}  y^{\frac{a}{2}}
\phi^{\frac{n+1}{n+a-1}} |\ana \phi|
|v| \dd     \nonumber         \\
 \geq
 c \left( \int_0^{+\infty} \int_{\R^n} |\phi^{\frac{2n+a}{n+a-1}} v|^{\frac{2(n+1)}{2n+a}}
\dd \right)^{\frac{2n+a}{2(n+1)}}   \ .
\eea
Next we will control the second term of the LHS by the first term of the LHS. To this end we consider
various cases. Suppose first that $s \in [\frac12,1)$ that is $a \in (-1,0]$ and $x \in \R^n_{+}$.
 Using the asymptotics of  Lemma  \ref{lem33} we get that
\[
 y^{\frac{a}{2}}  \phi^{\frac{n+1}{n+a-1}} |\ana \phi|
\sim
 \frac{ y^{\frac{a}{2}}}{(x_n^2+y^2)^{\frac{a(n+1)}{4(n+a-1)}+ \frac{a+2}{4}}} \ ,
\]
whereas,
\be\la{4.54}
 y^{\frac{a}{2}}  \phi^{\frac{2n+a}{n+a-1}} \sim
 \frac{ y^{\frac{a}{2}}  }{(x_n^2+y^2)^{\frac{a(2n+a)}{4(n+a-1)}}} \ .
\ee
We now apply  Lemma \ref{lem:abc} with the choice:
$A= \frac{a}{2}$, $B =0$ and $\xG = \frac{a(n+1)}{4(n+a-1)}+\frac{a+2}{4}$ taking into
account that
\[
A+B+2-2\xG = \frac{(2-a)(n-1)}{2(n+a-1)} >0  \ .
\]
Thus we get for some positive constant $c$ that \be\la{4.56}
\int_0^{+\infty} \int_{\R^n_{+}}   y^{\frac{a}{2}}
\phi^{\frac{2n+a}{n+a-1}}  |\ana v| \dd \geq c \int_0^{+\infty}
\int_{\R^n_{+}}  y^{\frac{a}{2}}  \phi^{\frac{n+1}{n+a-1}} |\ana
\phi| |v| \dd  \ . \ee We next consider the case $a \in (0,1)$, $x
\in \R^n_{+}$. In this case
\[
 y^{\frac{a}{2}}  \phi^{\frac{n+1}{n+a-1}} |\ana \phi|
\sim
 \frac{ y^{-\frac{a}{2}}}{(x_n^2+y^2)^{\frac{a(n+1)}{4(n+a-1)}+ \frac{2-a}{4}}} \ ,
\]
whereas,
\be\la{4.64}
 y^{\frac{a}{2}}  \phi^{\frac{2n+a}{n+a-1}} \sim
 \frac{ y^{\frac{a}{2}}  }{(x_n^2+y^2)^{\frac{a(2n+a)}{4(n+a-1)}}} \ .
\ee
We now use Lemma \ref{lem:r2} with the choice $A=\frac{a}{2}$, $B=0$  and $\xG= \frac12 +\frac{a(2n+a)}{4(n+a-1)}$
taking into account that
$\frac{x_n}{(x_n^2+y^2)^{\frac12}}<1$ and $A+B+2-2\xG = \frac{(2-a)(n-1)}{2(n+a-1)} >0 $.
 We then conclude that (\ref{4.56}) is valid for all $a \in (-1,1)$.

In a  similar manner for all $a \in (-1,1)$ and $x \in \R^n_{-}$ we get that
\[
 y^{\frac{a}{2}}  \phi^{\frac{n+1}{n+a-1}} |\ana \phi|
\sim
 \frac{ y^{-\frac{a}{2} +\frac{(1-a)(n+1)}{n+a-1}}}{(x_n^2+y^2)^{\frac{(2-a)(2n+a)}{4(n+a-1)}}} \ ,
\]
whereas,
\be\la{4.58}
 y^{\frac{a}{2}}  \phi^{\frac{2n+a}{n+a-1}} \sim
 \frac{ y^{\frac{a}{2}+\frac{(1-a)(2n+a)}{n+a-1}}   }{(x_n^2+y^2)^{\frac{(2-a)(2n+a)}{4(n+a-1)}}} \ .
\ee
This time we use Lemma \ref{lem:abc22} with $A=-\frac{a}{2}+\frac{(1-a)(n+1)}{n+a-1}$, $B=0$ and
$\xG  =\frac{(2-a)(2n+a)}{4(n+a-1)}$, noticing that
\[
A+B+2-2\xG = \frac{(2-a)(n-1)}{2(n+a-1)} >0  \ ,
\]
thus obtaining \be\la{4.62} \int_0^{+\infty} \int_{\R^n_{-}}
y^{\frac{a}{2}}  \phi^{\frac{2n+a}{n+a-1}}  |\ana v| \dd \geq c
\int_0^{+\infty} \int_{\R^n_{-}}  y^{\frac{a}{2}}
\phi^{\frac{n+1}{n+a-1}} |\ana \phi| |v| \dd  \ . \ee Combining
(\ref{4.56}) and (\ref{4.62}) we obtain the following $L^1$ Hardy
estimate on the whole $\R^n$: \be\la{4.63} \int_0^{+\infty}
\int_{\R^n}   y^{\frac{a}{2}}  \phi^{\frac{2n+a}{n+a-1}}  |\ana v|
\dd \geq c \int_0^{+\infty} \int_{\R^n}  y^{\frac{a}{2}}
\phi^{\frac{n+1}{n+a-1}} |\ana \phi| |v| \dd  \ . \ee Using this
 in (\ref{4.60}) we get that
\be\la{4.70} \int_0^{+\infty} \int_{\R^n}   y^{\frac{a}{2}}
\phi^{\frac{2n+a}{n+a-1}}  |\ana v| \dd
 \geq
 c \left( \int_0^{+\infty} \int_{\R^n} |\phi^{\frac{2n+a}{n+a-1}} v|^{\frac{2(n+1)}{2n+a}}
\dd \right)^{\frac{2n+a}{2(n+1)}}   \ .
\ee

To continue we next set in (\ref{4.70})  $v =
|w|^{\frac{2n+a}{n+a-1}}$ and apply Schwartz inequality in the LHS
to conclude after a simplification \be\la{4.72} \int_0^{+\infty}
\int_{\R^n}  y^{a} \phi^2 |\ana w|^2 \dd \geq c
\left(\int_0^{+\infty} \int_{\R^n} |\phi w|^{\frac{2(n+1)}{n+a-1}}
\dd \right)^{\frac{n+a-1}{n+1}} \ , \ee which is equivalent to
(\ref{4.48}). The result then follows.

\finedim

We are now ready to establish the  Proof of Theorem \ref{th16} part (ii).

\noindent {\em Proof of Theorem \ref{th16} part (ii):}  Again we
will use inequality (\ref{44hs}). This time the result will follow
once we will establish the following inequality: \be\la{4.74}
\int_0^{+\infty} \int_{\R^n}   y^a |\ana u - \frac{\ana
\phi}{\phi} u|^2 \dd \geq c \left( \int_{\R^n_{+}}
|u(x,0)|^{\frac{2n}{n+a-1}} dx \right)^{\frac{n+a-1}{n}} \ , \ee
with $\phi$   given by
\[
\phi(x,y)   =
 (y^2+x_n^2)^{-\frac{a}{4}} B(\frac{x_n}{y}),    ~~~~~~ y >0, ~~x_n \in \R \ ,
\]
and $B$ solves (\ref{ode4}), (\ref{bc2}).

 Our starting point  is  again the following  weighted  trace Sobolev
inequality,  see \cite{Maz}, Theorem 1, section 2.1.6,
 valid for functions $u \in C^{\infty}_{0}(\R^n \times \R)$ with $u(x,0)= 0$, $x \in \R^n_{-}$:
\[
\int_0^{+\infty} \int_{\R^n}  y^{\frac{a}{2}} |\ana u| \dd \geq c
\left( \int_{\R^n_{+}} |u(x,0)|^{\frac{2n}{2n+a}} dx
\right)^{\frac{2n+a}{2n}} \ .
\]
We set $u = \phi^{\frac{2n+a}{n+a-1}} v$ to obtain \bea\la{4.76}
\int_0^{+\infty} \int_{\R^n}   y^{\frac{a}{2}}
\phi^{\frac{2n+a}{n+a-1}}  |\ana v| \dd + \frac{2n+a}{n+a-1}
\int_0^{+\infty} \int_{\R^n}  y^{\frac{a}{2}}
\phi^{\frac{n+1}{n+a-1}} |\ana \phi|
|v| \dd     \nonumber         \\
 \geq
 c \left( \int_{\R^n_{+}} |\phi^{\frac{2n+a}{n+a-1}}(x,0) v(x,0)|^{\frac{2n}{2n+a}}
dx \right)^{\frac{2n+a}{2n}}   \ . \eea Combining this with
(\ref{4.63}) we get that \bea\la{4.78} \int_0^{+\infty}
\int_{\R^n} y^{\frac{a}{2}}  \phi^{\frac{2n+a}{n+a-1}}  |\ana v|
\dd
 \geq
 c \left( \int_{\R^n_{+}} |\phi^{\frac{2n+a}{n+a-1}}(x,0) v(x,0)|^{\frac{2n}{2n+a}}
dx \right)^{\frac{2n+a}{2n}}   \ .
\eea
We set  $v = |w|^{\frac{2n+a}{n+a-1}}$
 and apply Schwartz inequality in the
LHS to arrive at \bean \left( \int_0^{+\infty} \int_{\R^n}  y^{a}
\phi^2 |\ana w|^2 \dd \right)^{\frac12} \left(\int_0^{\infty}
\int_{\R^n} |\phi w|^{\frac{2(n+1)}{n+a-1}} \dd \right)^{\frac12}
\geq  c \left( \int_{\R^n_{+}} |(\phi w)(x,0) |^{\frac{2n}{n+a-1}}
dx \right)^{\frac{2n+a}{2n}}   \ . \eean We next use the Sobolev
inequality  (\ref{4.72}) to conclude after a simplification
\[
 \int_0^{+\infty} \int_{\R^n}  y^{a} \phi^2 |\ana w|^2 \dd
\geq c
\left( \int_{\R^n_{+}} |(\phi w)(x,0) |^{\frac{2n}{n+a-1}}
dx \right)^{\frac{n+a-1}{n}}   \ ,
\]
which is equivalent to (\ref{4.74}) and the result follows.

\finedim

\setcounter{equation}{0}
\section{The General Case,  Trace Hardy--Sobolev--Maz'ya I   \&  II}\la{sec5}

\subsection{Trace Hardy--Sobolev--Maz'ya I}
Here we will give the proof of Theorem \ref{th11} part (iii). We first establish the following Hardy--Sobolev--Maz'ya
where only the Hardy term appears in the trace term.

\begin{theorem}\la{th51}
  Let $\frac12 < s <1$, $n \geq 2$  and $\xO \subsetneqq \Ren$  be a uniformly Lipschitz domain  with finite inner
 radius that in addition satisfies
\be\la{5.2}
- \Delta d(x) \geq 0, ~~~~~~~~~~~~~ x \in \xO  \ .
\ee
 Then there exists a positive  constant $c$ such that   for all   $u \in C^{\infty}_{0}(\xO \times \R) $    there holds
\be\la{5.4}
 \int_{0}^{+\infty} \int_{\xO}  y^{1-2s} |\ana_{(x,y)} u(x,y)|^2 \dd  \geq  \bar{d}_{s}
 \int_{\xO} \frac{u^2(x,0)}{d^{2s}(x)}dx + c
 \left(\int_{0}^{+\infty}  \int_{\xO} | u(x,y)|^{\frac{2(n+1)}{n-2s}} dx dy \right)^{\frac{n-2s}{n+1}}  \ .
\ee
with
\be\la{5.6}
 \bar{d}_{s} :=  \frac{2 \Gamma \left(1-s \right) \Gamma^2 \left(\frac{3+2s}{4} \right)}
{\Gamma^2 \left(\frac{3-2s}{4} \right) \Gamma \left( s \right)} \ .
\ee
\end{theorem}

\noindent {\em Proof of Theorem \ref{th51}:} From the proof of Theorem \ref{th11} we recall the inequality
(\ref{id2}), that is
\bea\la{5.8}
 \dino y^a |\ana u|^2 \dd     \geq  &   &
\bar{d}_s  \int_{\xO}  \frac{ u^2(x,0)}{d^{1-a}(x)} dx
+ \dino y^a |\ana u - \frac{\ana \phi}{\phi} u|^2 \dd  \nonumber \\
&  - &
 \dino \frac{{\rm div}(y^a \ana \phi)}{\phi} u^2 \dd  \ ,
\eea
where
 $\phi$ is  given by
\be\la{5.9}
\phi(x,y)= d^{-\frac{a}{2}}(x) A \left( \frac{y}{d} \right), ~~~~~~~~
 y >0,~~~~~ x \in \xO  \ ,
\ee
and  $A$ solves (\ref{ode1}), (\ref{2bc}).

The result will follow after  establishing  the following inequality:
\be\la{5.10}
\dino y^a |\ana u - \frac{\ana \phi}{\phi} u|^2 \dd
- \dino \frac{{\rm div}(y^a \ana \phi)}{\phi} u^2 \dd \geq
 c
 \left(\int_{0}^{+\infty}  \int_{\xO} | u(x,y)|^{\frac{2(n+1)}{n+a-1}} dx dy \right)^{\frac{n+a-1}{n+1}}  \ .
\ee

To this end we start with  the inequality, see \cite{Maz}, Theorem 1, section 2.1.6,
\[
\int_0^{+\infty} \int_{\xO}  y^{\frac{a}{2}} |\ana u| \dd \geq c
\left(\int_0^{+\infty} \int_{\xO} |u(x,y)|^{\frac{2(n+1)}{2n+a}}
\dd \right)^{\frac{2n+a}{2(n+1)}} \ ,~~~~~~~ u \in
C^{\infty}_{0}(\xO \times \R) \ ,
\]
with the choice $u = \phi^{\frac{2n+a}{n+a-1}} v$. Hence we obtain
\bea\la{5.12} \int_0^{+\infty} \int_{\xO}   y^{\frac{a}{2}}
\phi^{\frac{2n+a}{n+a-1}}  |\ana v| \dd + \frac{2n+a}{n+a-1}
\int_0^{+\infty} \int_{\xO}  y^{\frac{a}{2}}
\phi^{\frac{n+1}{n+a-1}} |\ana \phi|
|v| \dd     \nonumber         \\
 \geq
 c \left( \int_0^{+\infty} \int_{\xO} |\phi^{\frac{2n+a}{n+a-1}} v|^{\frac{2(n+1)}{2n+a}}
\dd \right)^{\frac{2n+a}{2(n+1)}}   \ .
\eea
Next we will control the second term of the LHS using Lemma \ref{lem:gen}. To this end we recall that
for $a \in (-1,0)$ we have the following asymptotics from Lemma  \ref{lemasf}:
\[
 y^{\frac{a}{2}}  \phi^{\frac{n+1}{n+a-1}} |\ana \phi|
\sim
 \frac{ y^{\frac{a}{2}}  d^{\frac{n+1}{n+a-1}}}{(d^2+y^2)^{\frac{(2+a)(2n+a)}{4(n+a-1)}}} \ ,
\]
whereas,
\be\la{5.14}
 y^{\frac{a}{2}}  \phi^{\frac{2n+a}{n+a-1}} \sim
 \frac{ y^{\frac{a}{2}}  d^{\frac{2n+a}{n+a-1}}}{(d^2+y^2)^{\frac{(2+a)(2n+a)}{4(n+a-1)}}} \ .
\ee
We then use Lemma \ref{lem:gen} with the choice $A= \frac{a}{2}$,  $B = \frac{n+1}{n+a-1}$ and
 $\xG = \frac{(2+a)(2n+a)}{4(n+a-1)}$  taking into
account that
\[
A+B+2-2\xG = \frac{(2-a)(n-1)}{2(n+a-1)} >0  \ ,
\]
to obtain the estimate
\bean
 \int_0^{+\infty} \int_{\xO}
 \frac{ y^{\frac{a}{2}}  d^{\frac{n+1}{n+a-1}}}{(d^2+y^2)^{\frac{(2+a)(2n+a)}{4(n+a-1)}}}
|v| \dd  \leq  C_1   \int_0^{+\infty} \int_{\xO}
 \frac{ y^{\frac{a}{2}}  d^{\frac{2n+a}{n+a-1}}}{(d^2+y^2)^{\frac{(2+a)(2n+a)}{4(n+a-1)}}}  | \nabla v| \dd   \\
+ C_2  \int_0^{+\infty} \int_{\xO}
 \frac{ y^{\frac{a}{2}}  d^{\frac{2n+a}{n+a-1}}}{(d^2+y^2)^{\frac{(2+a)(2n+a)}{4(n+a-1)}}}  | v| \dd   \ .
\eean From this and (\ref{5.12}) we have that
\[
\int_0^{+\infty} \int_{\xO}   y^{\frac{a}{2}}
\phi^{\frac{2n+a}{n+a-1}}  |\ana v| \dd  + \int_0^{+\infty}
\int_{\xO}   y^{\frac{a}{2}}  \phi^{\frac{2n+a}{n+a-1}}  | v| \dd
\geq
 c \left( \int_0^{+\infty} \int_{\xO} |\phi^{\frac{2n+a}{n+a-1}} v|^{\frac{2(n+1)}{2n+a}}
\dd \right)^{\frac{2n+a}{2(n+1)}}
\]
To continue we next set  $v = |w|^{\frac{2n+a}{n+a-1}}$ and apply
Schwartz inequality in the LHS. After a simplification  we arrive
at:
 \be\la{5.13} \dino y^a \phi^2 |\nabla w|^2 \dd + \dino y^a
\phi^2  w^2  \dd \geq C
 \left( \dino |\phi w|^{\frac{2(n+1)}{n+a-1}} \right)^{\frac{n+a-1}{n+1}}
\ee
To conclude the proof of the Theorem we need the following estimate:
\be\la{5.16}
c  \dino y^a \phi^2  w^2  \dd  \leq \dino y^a \phi^2 |\nabla w|^2 \dd  -
 \dino {\rm div}(y^a \nabla \phi) \phi w^2 \dd  \ .
\ee
It is here that we will  use the fact that the domain $\xO$ has finite inner radius.
Using  Lemma \ref{lem:l2log} with $A=a$, $B=0$ we obtain that
\[
 c \dino \frac{ y^a X^2\left(\frac{d}{R_{in}}\right)}{(d^2+y^2)^{\frac{2+a}{2}}} w^2 \dd \leq
 \dino \frac{ y^a d^2}{(d^2+y^2)^{\frac{2+a}{2}}} | \nabla w|^2 \dd  -
 \dino \frac{ y^a d (\Delta d) X\left(\frac{d}{R_{in}}\right)  }{(d^2+y^2)^{\frac{2+a}{2}}}  w^2 \dd \ ,
\]
which implies
\[
c  \dino \frac{ y^a d^2}{(d^2+y^2)^{\frac{2+a}{2}}} w^2 \dd \leq
 \dino \frac{ y^a d^2}{(d^2+y^2)^{\frac{2+a}{2}}} | \nabla w|^2 \dd  -
 \dino \frac{ y^a d (\Delta d)}{(d^2+y^2)^{\frac{2+a}{2}}}  w^2 \dd \ .
\]
Taking into account the asymptotics of $\phi$
this is equivalent to (\ref{5.16}). We omit further details.

\finedim

We are now ready to prove Theorem \ref{th11} part (iii).

\noindent {\em Proof of Theorem \ref{th11} part (iii):} Again we will use (\ref{5.8}).
 The result then will follow  once  we  establish:
\be\la{5.18}
\dino y^a |\ana u - \frac{\ana \phi}{\phi} u|^2 \dd
- \dino \frac{{\rm div}(y^a \ana \phi)}{\phi} u^2 \dd \geq
 c
 \left(  \int_{\xO} | u(x,0)|^{\frac{2n}{n+a-1}} dx  \right)^{\frac{n+a-1}{n}}  \ .
\ee
where $\phi$ is as in (\ref{5.9}).
To this end we start with  the inequality, see \cite{Maz}, Theorem 1, section 2.1.6,
\[
\int_0^{+\infty} \int_{\xO}  y^{\frac{a}{2}} |\ana u| \dd \geq c
\left( \int_{\xO} |u(x,0)|^{\frac{2n}{2n+a}} dx
\right)^{\frac{2n+a}{2n}} \ ,~~~~~~~ u \in C^{\infty}_{0}(\xO
\times \R) \ ,
\]
with the choice $u = \phi^{\frac{2n+a}{n+a-1}} v$. Hence we obtain
\bea\la{5.19}
 \int_0^{+\infty} \int_{\xO}   y^{\frac{a}{2}}
\phi^{\frac{2n+a}{n+a-1}}  |\ana v| \dd + \frac{2n+a}{n+a-1}
\int_0^{+\infty} \int_{\xO}  y^{\frac{a}{2}}
\phi^{\frac{n+1}{n+a-1}} |\ana \phi|
|v| \dd     \nonumber         \\
 \geq
 c \left(  \int_{\xO} |\phi^{\frac{2n+a}{n+a-1}} v|^{\frac{2n}{2n+a}}
dx \right)^{\frac{2n+a}{2n}}   \ .
\eea
Next we will control the second term of the LHS exactly as we did in the proof of Theorem \ref{th51}, to arrive at
\[
\int_0^{+\infty} \int_{\xO}   y^{\frac{a}{2}}
\phi^{\frac{2n+a}{n+a-1}}  |\ana v| \dd  + \int_0^{+\infty}
\int_{\xO}   y^{\frac{a}{2}}  \phi^{\frac{2n+a}{n+a-1}}  | v| \dd
\geq
 c \left(  \int_{\xO} |\phi^{\frac{2n+a}{n+a-1}} v(x,0)|^{\frac{2n}{2n+a}}
dx \right)^{\frac{2n+a}{2n}}  \ .
\]
To continue we next set  $v = |w|^{\frac{2n+a}{n+a-1}}$ and apply Schwartz inequality in the
LHS to get after elementary manipulations that
\bea\la{5.20}
 \left( \dino |\phi w|^{\frac{2(n+1)}{n+a-1}} \dd  \right)
 \left[  \dino y^a \phi^2 |\nabla w|^2 \dd + \dino y^a \phi^2  w^2  \dd \right]  \nonumber   \\
  \geq C
 \left(  \int_{\xO}  |\phi w(x,0)|^{\frac{2n}{n+a-1}} dx \right)^{\frac{2n+a}{n}}  \ .
\eea
At this point we use Theorem \ref{th51} and inequality (\ref{5.16}) to conclude the result. We
omit further details.

\finedim

\subsection{Trace Hardy--Sobolev--Maz'ya II}
Here we will give the proof of Theorem \ref{th14} part (iii). We first establish the following Hardy--Sobolev--Maz'ya
where only the Hardy term appears in the trace term.

\begin{theorem}\la{th52}
  Let $\frac12 < s <1$, $n \geq 2$  and $\xO \subsetneqq \Ren$  be a uniformly Lipschitz and convex
  domain with finite inner
 radius.
 Then, there exists a positive  constant $c$ such that   for all   $u \in C^{\infty}_{0}(\R^n \times \R) $ with
$u(x,0)=0$ for $x \in \CC \xO$
   there holds
\be\la{5.24}
 \int_{0}^{+\infty} \int_{\R^n}  y^{1-2s} |\ana_{(x,y)} u(x,y)|^2 \dd  \geq  \bar{k}_{s}
 \int_{\xO} \frac{u^2(x,0)}{d^{2s}(x)}dx + c
 \left(\int_{0}^{+\infty}  \int_{\R^n} | u(x,y)|^{\frac{2(n+1)}{n-2s}} dx dy \right)^{\frac{n-2s}{n+1}}  \ .
\ee
with
\be\la{5.26}
 \bar{k}_{s} :=  \frac{ 2^{1-2 s}  \Gamma^2(s +\frac12) \Gamma(1-s)}{\pi \Gamma(s)}  \ .
\ee

\end{theorem}

\noindent {\em Proof of Theorem \ref{th52}:} From the proof of Theorem \ref{th14} we recall the inequality
(\ref{id2b}), that is
\bea\la{5.28}
 \dinr y^a |\ana u|^2 \dd     \geq  &   &
\bar{k}_s  \int_{\xO}  \frac{ u^2(x,0)}{d^{1-a}(x)} dx
+ \dinr y^a |\ana u - \frac{\ana \phi}{\phi} u|^2 \dd  \nonumber \\
&  - &
 \dinr \frac{{\rm div}(y^a \ana \phi)}{\phi} u^2 \dd  \ ,
\eea
where
 $\phi$ is  given by
\be\la{5.29}
\phi(x,y)   = \left\{ \begin{array}{ll}
 (y^2+d^2)^{-\frac{a}{4}} B(\frac{d}{y}),  &  ~~~~~~~  x \in \xO,~~~ y >0     \\
  (y^2+d^2)^{-\frac{a}{4}} B(-\frac{d}{y}), & ~~~~~~~ x \in \CC \xO,~~ y >0  \ ,  \\
\end{array} \right.
\ee
and $B$ is the solution of the boundary value problem (\ref{ode4}) and (\ref{bc2}).
The result will follow after  establishing  the following inequality:
\be\la{5.30}
\dinr y^a |\ana u - \frac{\ana \phi}{\phi} u|^2 \dd
- \dinr \frac{{\rm div}(y^a \ana \phi)}{\phi} u^2 \dd \geq
 c
 \left(\int_{0}^{+\infty}  \int_{\R^n} | u(x,y)|^{\frac{2(n+1)}{n+a-1}} dx dy \right)^{\frac{n+a-1}{n+1}}  \ .
\ee
To this end we start with  the inequality, see \cite{Maz}, Theorem 1, section 2.1.6,
\[
\int_0^{+\infty} \int_{\R^n}  y^{\frac{a}{2}} |\ana u| \dd \geq c
\left(\int_0^{+\infty} \int_{\R^n} |u(x,y)|^{\frac{2(n+1)}{2n+a}}
\dd \right)^{\frac{2n+a}{2(n+1)}} \ ,~~~~~~~ u \in
C^{\infty}_{0}(\R^n \times \R) \ ,
\]
with the choice $u = \phi^{\frac{2n+a}{n+a-1}} v$.
 Hence we obtain
\bea\la{5.32} \int_0^{+\infty} \int_{\R^n}   y^{\frac{a}{2}}
\phi^{\frac{2n+a}{n+a-1}}  |\ana v| \dd + \frac{2n+a}{n+a-1}
\int_0^{+\infty} \int_{\R^n}  y^{\frac{a}{2}}
\phi^{\frac{n+1}{n+a-1}} |\ana \phi|
|v| \dd     \nonumber         \\
 \geq
 c \left( \int_0^{+\infty} \int_{\R^n} |\phi^{\frac{2n+a}{n+a-1}} v|^{\frac{2(n+1)}{2n+a}}
\dd \right)^{\frac{2n+a}{2(n+1)}}   \ .
\eea
Again we want to control  the second term of the LHS. This time we split the integral into the integral over
$\xO$ and the integral over  $\CC \xO$. Concerning the integral over  $\CC \xO$ we use the asymptotics of $\phi$
as given by Lemma \ref{lem33} for $a \in (-1,0)$  to get that
\[
 y^{\frac{a}{2}}  \phi^{\frac{n+1}{n+a-1}} |\ana \phi|
\sim
 \frac{ y^{-\frac{a}{2} +\frac{(1-a)(n+1)}{n+a-1}}}{(d^2+y^2)^{\frac{(2-a)(2n+a)}{4(n+a-1)}}} \ ,
\]
whereas,
\[
 y^{\frac{a}{2}}  \phi^{\frac{2n+a}{n+a-1}} \sim
 \frac{ y^{\frac{a}{2}+\frac{(1-a)(2n+a)}{n+a-1}}   }{(d^2+y^2)^{\frac{(2-a)(2n+a)}{4(n+a-1)}}} \ .
\]
This time we use Lemma \ref{lem:abc22} with $A=-\frac{a}{2}+\frac{(1-a)(n+1)}{n+a-1}$, $B=0$ and
$\xG  =\frac{(2-a)(2n+a)}{4(n+a-1)}$, noticing that
\[
A+B+2-2\xG = \frac{(2-a)(n-1)}{2(n+a-1)} >0  \ ,
\]
thus obtaining \be\la{5.34} \int_0^{+\infty} \int_{\CC \xO }
y^{\frac{a}{2}}  \phi^{\frac{2n+a}{n+a-1}}  |\ana v| \dd \geq c
\int_0^{+\infty} \int_{\CC \xO}  y^{\frac{a}{2}}
\phi^{\frac{n+1}{n+a-1}} |\ana \phi| |v| \dd  \ , \ee where we
also used the convexity of $\xO$.

On the other hand in $\xO$ the   asymptotics of $\phi$ are also given  by Lemma \ref{lem33} as follows:
\[
 y^{\frac{a}{2}}  \phi^{\frac{n+1}{n+a-1}} |\ana \phi|
\sim
 \frac{ y^{\frac{a}{2}}}{(d^2+y^2)^{\frac{a(2n+a)}{4(n+a-1)}+ \frac{1}{2}}} \ ,
\]
whereas,
\[
 y^{\frac{a}{2}}  \phi^{\frac{2n+a}{n+a-1}} \sim
 \frac{ y^{\frac{a}{2}}  }{(d^2+y^2)^{\frac{a(2n+a)}{4(n+a-1)}}} \ .
\]
We next use  Lemma \ref{lem:gen} with the choice $A= \frac{a}{2}$,  $B=0$ and
 $\xG = \frac{a(2n+a)}{4(n+a-1)}+ \frac{1}{2}$
taking into
account that
\[
A+B+2-2\xG = \frac{(2-a)(n-1)}{2(n+a-1)} >0  \ ,
\]
to obtain the estimate \bea\la{5.36} & &  \int_0^{+\infty}
\int_{\xO}
 \frac{ y^{\frac{a}{2}} }{(d^2+y^2)^{\frac{a(2n+a)}{4(n+a-1)}+ \frac{1}{2}}}
|v| \dd    \nonumber   \\
   &  \leq &   C_1   \int_0^{+\infty} \int_{\xO}
 \frac{ y^{\frac{a}{2}}  d}{(d^2+y^2)^{\frac{a(2n+a)}{4(n+a-1)}+ \frac{1}{2}} }  | \nabla v| \dd
+ C_2  \int_0^{+\infty} \int_{\xO}
 \frac{ y^{\frac{a}{2}}  d}{(d^2+y^2)^{\frac{a(2n+a)}{4(n+a-1)}+ \frac{1}{2}} }  | v| \dd \nonumber   \\
&  \leq &   C_1   \int_0^{+\infty} \int_{\xO}
 \frac{ y^{\frac{a}{2}}  }{(d^2+y^2)^{\frac{a(2n+a)}{4(n+a-1)}} }  | \nabla v| \dd
+ C_2  \int_0^{+\infty} \int_{\xO}
 \frac{ y^{\frac{a}{2}}  d}{(d^2+y^2)^{\frac{a(2n+a)}{4(n+a-1)}+ \frac{1}{2}} }  | v| \dd  \ .\nonumber
\eea
Equivalently, this can be written as
\bea\la{5.38}
&  & C \int_0^{+\infty} \int_{\xO}  y^{\frac{a}{2}} \phi^{\frac{n+1}{n+a-1}} |\nabla \phi| |v| \dd  \nonumber \\
& \leq &  \int_0^{+\infty} \int_{\xO}  y^{\frac{a}{2}}
\phi^{\frac{2n+a}{n+a-1}} |\nabla v| \dd + \int_0^{+\infty}
\int_{\xO}  y^{\frac{a}{2}} \phi^{\frac{2n+a}{n+a-1}}
\frac{d}{(d^2+y^2)^{\frac12}} |v| \dd  \ . \eea Using (\ref{5.34})
and (\ref{5.38})  in (\ref{5.32}) we arrive at \bea\la{5.40}
\int_0^{+\infty} \int_{\R^n}   y^{\frac{a}{2}}
\phi^{\frac{2n+a}{n+a-1}}  |\ana v| \dd +\int_0^{+\infty}
\int_{\xO}  y^{\frac{a}{2}}  \frac{d}{(d^2+y^2)^{\frac12}}
\phi^{\frac{2n+a}{n+a-1}}
|v| \dd     \nonumber         \\
 \geq
 c \left( \int_0^{+\infty} \int_{\R^n} |\phi^{\frac{2n+a}{n+a-1}} v|^{\frac{2(n+1)}{2n+a}}
\dd \right)^{\frac{2n+a}{2(n+1)}}   \ .
\eea
To continue we next set  $v = |w|^{\frac{2n+a}{n+a-1}}$ and apply Schwartz inequality in the
LHS. After a simplification  we arrive at:
\be\la{5.42}
\dinr y^a \phi^2 |\nabla w|^2 \dd + \dino  \frac{ y^a   d^2 \phi^2   }{d^2+y^2}   w^2  \dd \geq c
 \left( \dinr |\phi w|^{\frac{2(n+1)}{n+a-1}} \right)^{\frac{n+a-1}{n+1}}
\ee
To conclude the proof of the Theorem it is enough  to obtain the following estimate:
\be\la{5.44}
c   \dino  \frac{ y^a   d^2 \phi^2   }{d^2+y^2}   w^2  \dd
 \leq \dino y^a \phi^2 |\nabla w|^2 \dd  -
 \dino {\rm div}(y^a \nabla \phi) \phi w^2 \dd  \ .
\ee

It is here that we will  use the fact that the domain $\xO$ has finite inner radius.
Using  Lemma \ref{lem:l2log} with $A=a$, $B=0$ we obtain that
\[
 c \dino \frac{ y^a X^2\left(\frac{d}{R_{in}}\right)}{(d^2+y^2)^{\frac{2+a}{2}}} w^2 \dd \leq
 \dino \frac{ y^a d^2  }{(d^2+y^2)^{\frac{2+a}{2}}} | \nabla w|^2 \dd  -
 \dino \frac{ y^a d (\Delta d) X\left(\frac{d}{R_{in}}\right)  }{(d^2+y^2)^{\frac{2+a}{2}}}  w^2 \dd \ ,
\]
which implies
\[
c  \dino \frac{ y^a d^2}{(d^2+y^2)^{\frac{2+a}{2}}} w^2 \dd \leq
 \dino \frac{ y^a }{(d^2+y^2)^{\frac{a}{2}}} | \nabla w|^2 \dd  -
 \dino \frac{ y^a  (\Delta d)}{(d^2+y^2)^{\frac{1+a}{2}}}  w^2 \dd \ .
\]
Taking into account the asymptotics of $\phi$
this is equivalent to (\ref{5.44}). We omit further details.

\finedim

We are now ready to prove Theorem \ref{th14} part (iii).

\noindent {\em Proof of Theorem \ref{th14} part (iii):} Again we will use (\ref{5.28}).
 The result then will follow  once  we  establish:
\be\la{5.46}
\dinr y^a |\ana u - \frac{\ana \phi}{\phi} u|^2 \dd
- \dinr \frac{{\rm div}(y^a \ana \phi)}{\phi} u^2 \dd \geq
 c
 \left(  \int_{\xO} | u(x,0)|^{\frac{2n}{n+a-1}} dx  \right)^{\frac{n+a-1}{n}}  \ .
\ee
where $\phi$ is as in (\ref{5.29}).
To this end we start again  with  the inequality,
\[
\int_0^{+\infty} \int_{\R^n}  y^{\frac{a}{2}} |\ana u| \dd \geq c
\left( \int_{\xO} |u(x,0)|^{\frac{2n}{2n+a}} dx
\right)^{\frac{2n+a}{2n}} \ ,
\]
 valid for   $u \in C^{\infty}_{0}(\R^n \times \R)$ with
$u(x,0)=0,~~ x \in \CC \xO$. We apply this to  $u =
\phi^{\frac{2n+a}{n+a-1}} v$. Hence we obtain \bea\la{5.48}
\int_0^{+\infty} \int_{\R^n}   y^{\frac{a}{2}}
\phi^{\frac{2n+a}{n+a-1}}  |\ana v| \dd + \frac{2n+a}{n+a-1}
\int_0^{+\infty} \int_{\R^n}  y^{\frac{a}{2}}
\phi^{\frac{n+1}{n+a-1}} |\ana \phi|
|v| \dd     \nonumber         \\
 \geq
 c \left(  \int_{\xO} |\phi^{\frac{2n+a}{n+a-1}} v|^{\frac{2n}{2n+a}}
dx \right)^{\frac{2n+a}{2n}}   \ .
\eea
Next we will control the second term of the LHS exactly as we did in the proof of Theorem \ref{th52}, to arrive at
\[
\int_0^{+\infty} \int_{\R^n}   y^{\frac{a}{2}}
\phi^{\frac{2n+a}{n+a-1}}  |\ana v| \dd  + \int_0^{+\infty}
\int_{\xO}  \frac{  y^{\frac{a}{2}} d}{(d^2+y^2)^{\frac12}}
\phi^{\frac{2n+a}{n+a-1}}  | v| \dd \geq
 c \left(  \int_{\xO} |\phi^{\frac{2n+a}{n+a-1}} v(x,0)|^{\frac{2n}{2n+a}}
dx \right)^{\frac{2n+a}{2n}}  \ .
\]
To continue we next set  $v = |w|^{\frac{2n+a}{n+a-1}}$ and apply
Schwartz inequality in the LHS to get after elementary
manipulations that \bea\la{5.201}
 \left( \dino |\phi w|^{\frac{2(n+1)}{n+a-1}} \dd  \right)
 \left[  \dino y^a \phi^2 |\nabla w|^2 \dd + \dino \frac{y^a d^2}{d^2+y^2} \phi^2  w^2  \dd \right]  \nonumber   \\
  \geq C
 \left(  \int_{\xO}  |\phi w(x,0)|^{\frac{2n}{n+a-1}} dx \right)^{\frac{2n+a}{n}}  \ .
\eea
At this point we use Theorem \ref{th52} and inequality (\ref{5.44}) to conclude the result. We
omit further details.

\finedim

\setcounter{equation}{0}
\section{The Fractional Laplacians}
In this section we will apply the previous results to establish   the proofs of Theorems \ref{th13},  \ref{th15}
as well as of part (iii) of Theorem \ref{th16}.

\noindent {\em Proof of Theorem \ref{th13}:} Part (i) and (iii) follow from part (i) and (iii) of  Theorem \ref{th11}
taking into account the relation between the energy of the extended problem and the corresponding one of the
fractional Laplacian, see subsection  \ref{subsec81} and in particular relation   (\ref{ap1.4}).

We next prove part (ii).
  We will use the optimality of the constant $\bar{d}_s$ of
 Theorem \ref{th11}, that is
for each $\xe>0$ there exists a $u_{\xe} \in C^{\infty}_{0}(\xO \times \R)$ such that
\[
 \bar{d}_{s} + \xe  \geq
 \frac{\dino y^{1-2s} |\ana u_{\xe}|^2 \dd}{\int_{\xO} \frac{u^2_{\xe}(x,0)}{d^{2s}(x)} dx}  \ ,
\]
and let $f_{\xe}(x)=u_{\xe}(x,0)$. We will show that for some
positive constant $c$, \be\la{5.400}
 d_{s} + c \xe  \geq
 \frac{ ((-\Delta)^{s} f_{\xe} , f_{\xe} )_{\xO}}{\int_{\xO} \frac{f^2_{\xe}(x)}{d^{2s}(x)} dx}  \ .
\ee
To this end let $\hat{u}_{\xe}$ be the solution to the extended problem
\bean
{\rm div} (y^{1-2s} \ana \hat{u}_{\xe}(x,y)) & = & 0,  ~~~~~~~~  {\rm in} ~~~~~\xO \times (0, \infty) \ , \\
\hat{u}_{\xe}(x,y) & = & 0, ~~~~~~~~~x \in \partial \xO \times (0, \infty)  \ ,  \\
\hat{u}_{\xe}(x,0) & = & f_{\xe}(x)    \ .
\eean
The solution  $\hat{u}_{\xe}$ minimizes the energy and therefore
\[
\dino y^{1-2s} |\ana \hat{u}_{\xe}|^2 \dd  \leq \dino y^{1-2s} |\ana u_{\xe}|^2 \dd  \ .
\]
On the other hand using (\ref{ap1.4})  we have
\[
\dino y^{1-2s} |\ana \hat{u}_{\xe}|^2 \dd  =
\frac{2^{1-2s} \xG(1-s)}{ \xG(s)}    ((-\Delta)^{s} f_{\xe} , f_{\xe} )_{\xO} \ ,
\]
and (\ref{5.400}) follows easily  with $c = \frac{
\xG(s)}{2^{1-2s} \xG(1-s)}$.

\finedim

We next give the proof of Theorem  \ref{th15}

\noindent {\em Proof of Theorem \ref{th15}:} Part (i) and (iii) follow from part (i) and (iii) of  Theorem \ref{th14}
taking into account the relation between the energy of the extended problem and the corresponding one of the
fractional Laplacian, see subsection  \ref{subsec82} and in particular relations   (\ref{81})--(\ref{82}).

The proof of part (ii) is quite similar to the proof of part (ii) of Theorem \ref{th13}, the only difference
being that the extension problem is now  on the whole $\R^n$.  We  omit  the details.

\finedim

Finally estimate (\ref{1.C2})
  of part (iii) of Theorem  \ref{th16} follows at once from part (ii) of Theorem  \ref{th16}
and (\ref{81}).  Concerning estimate (\ref{1.C3}), it  follows from (\ref{1.C2}) taking into
account that for $x \in \R^n_{+}$,
\[
\int_{\R^n_{-}} \frac{ d \xi}{|x-\xi|^{n+2s}} =
\frac{ \pi^{\frac{n-1}{2}} \Gamma\left( \frac{1+2s}{2} \right)}{2s  \Gamma\left( \frac{n+2s}{2} \right)} \frac{1}{x_n^{2s}} \ ,
\]
see, e.g., \cite{BBC}.

\setcounter{equation}{0}
\section{Appendix}\la{sec8}

\subsection{Spectral Fractional  Laplacian}\la{subsec81}

Let $\xO \subset \Ren$ be a bounded  domain, and let  $\xl_i$ and  $\phi_i$ be the Dirichlet eigenvalues and
 eigenfunctions of the Laplacian,
i.e. $-\Delta \phi_i = \xl_i \phi_i$ in $\xO$, with $\phi_i = 0$ on $\partial \xO$, normalized so that
 $\int_{\xO}  \phi_i^2 dx = 1$. Then,  for $f(x) = \sum c_i \phi_i(x)$ we  define
\be\la{frc1}
(-\Delta)^{s} f =  \sum_{i=1}^{\infty} c_i \xl_i^s  \phi_i,~~~~~~~~~~~~0<s<1.
\ee
We also have
\be\la{fr2}
((-\Delta)^{s} f , f )_{\xO} = \int_{\xO} f~ (-\Delta)^{s} f dx =  \sum_{i=1}^{\infty}  c_i^2 \xl_i^s.
\ee

To the function $f(x)$ we associate the ``extended'' function $u(x,y)$, $x \in \xO$, $y >0$, given by
\[
u(x,y) = \sum_{i=1}^{+\infty} c_i \phi_i(x) T(y \sqrt{\xl_i}),
\]
where $T(t)$ is the energetic  solution of the ODE:
\be\la{app.14}
(t^{1-2s} T'(t))'-t^{1-2s} T(t) = 0,~~~~~{\rm or}~~~~ T''+\frac{1-2s}{t}T'-T=0,~~~~~~         t \geq 0.
\ee
The solution of this can be taken from \cite{AS}, Section 9.6 and is given by
\be\la{app.16}
T(t) = \frac{2^{1-s}}{\xG(s)} t^s K_s(t),
\ee
where $K_s(t)$ denotes the modified Bessel function of second kind. The constant factor is chosen in such a way
that $T(0)=1$. As a consequence we also have $u(x,0)=f(x)$.

An easy calculation shows that ${\rm div}(y^{1-2s} \ana (\phi_i(x) T(y \sqrt{\xl_i}))=0$ from which it follows that
${\rm div} (y^{1-2s} \ana u)=0$. An integration by parts then shows that
\bea\la{ap1.4}
\dino y^{1-2s} |\ana u|^2 \dd &  = &  \lim_{\tau \rft}  \tau^{1-2s} \int_{\xO} u(x,\tau)u_{y}(x,\tau) dx
 -\lim_{\tau  \ra 0} \tau^{1-2s} \int_{\xO} u(x,\tau)u_{y}(x,\tau) dx   \nonumber  \\
 &  = & \left[  \lim_{t \rft} t^{1-2s}T(t)T'(t) -
\lim_{t \ra 0} t^{1-2s}T(t)T'(t) \right] \sum_{i=1}^{\infty}  \xl_i^{s} c_i^2   \nonumber   \\
 &  = & \frac{2^{1-2s} \xG(1-s)}{ \xG(s)}    ((-\Delta)^{s} f , f )_{\xO} \ .
\eea
Where we used (\ref{fr2}) and the fact that
\be\la{app.18}
  \lim_{t \rft} t^{1-2s}T(t)T'(t) -
\lim_{t \ra 0} t^{1-2s}T(t)T'(t) = \frac{2^{1-2s} \xG(1-s)}{ \xG(s)} \ .
\ee
To prove the above relation we show that
\[
 \lim_{t \rft} t^{1-2s}T(t)T'(t) =0, ~~~~~~~~~~~-\lim_{t \ra 0} t^{1-2s}T(t)T'(t) = \frac{2^{1-2s} \xG(1-s)}{ \xG(s)}.
\]
These two relations are a direct  consequence of (\ref{app.16}) and  the
 following properties of $K_s(t)$ :
\bean
K_s(t) &  \sim  & \frac{\xG(s)}{2^{1-s}}t^{-s},~~~ t \ra 0,~~~~~~ K_s(t) \sim \sqrt{\frac{\pi}{2t}} e^{-t}, ~~~~~~t \rft, \\
\frac{d}{dt}(t^s K_s(t))&  =  & -t^s K_{s-1}(t), ~~~~~~~~~~~~~~~~K_s(t)=K_{-s}(t) \ .
\eean

\subsection{Dirichlet Fractional Laplacian}\la{subsec82}

Let $u(x,y)$ be the extended function as defined in  (\ref{ext1})--(\ref{ext2}).
In this subsection we will show  the following two relations connecting the energy of the extended problem
and the energy of the  Dirichlet fractional Laplacian:

\bea \la{81}
 \int^{+\infty}_{0} \int_{\R^n} y^{1-2s}  |\nabla u|^2  dx dy &   =  &
\frac{s\Gamma\left(\frac{n+2s}{2}\right)}{\pi^{\frac{n}{2}} \Gamma(s)}
\int_{\R^n} \int_{\R^n}
\frac{|f(x)-f(\xi)|^2}{|x-\xi|^{n+2s}} dx d\xi     \\
\dinr y^{1-2s} |\ana u|^2 \dd &  = &
\frac{2^{1-2s} \xG(1-s)}{ \xG(s)}    ((-\Delta)^{s} f , f )_{\Ren}   \la{82}
 \ .
\eea

We will use the Fourier transform in the $x$-variables:
\[
\hat{u}(\eta, y) = (2 \pi)^{-\frac{n}{2}} \int_{\Ren} e^{-i x \cdot \eta} u(x,y) dx.
\]
The  equation ${\rm div}(y^{1-2s} \nabla u(x,y))=0$ or equivalently
$\Delta_x u+ u_{yy}+\frac{a}{y} u_y=0 $  with $u(x,0)=f(x)$, reads as follows when taking
the Fourier transform
\[
- |\eta|^2 \hat u + (\hat
u)_{yy}+\frac{1-2s}{y} (\hat u)_{y}=0, ~~~~~~~~\hat u(\eta, 0)= \hat f(\eta),
\]
 and it is satisfied by $\hat
u(\eta, y)=\hat f(\eta) T(|\eta|y)$, where $T$ satisfies (\ref{app.14}) and is given by (\ref{app.16}).

Concerning the energies  we have: \bean \int^{+\infty}_{0}
y^{1-2s} \int_{\R^n} |\nabla u|^2 dx dy & = & \int^{+\infty}_{0}
y^{1-2s} \int_{\R^n} \left( |\eta|^2  |\hat { u}|^2 +
|\hat{u}_y|^2 \right) d\eta dy  \\
&  =  & \int^{+\infty}_{0} y^{1-2s} \int_{\R^n} \left\{ |\eta|^2
|\hat f|^2 [T^2(|\eta|y) + T^{'2}(|\eta|y)]
\right\} d\eta dy  \\
 &  =  &  \left(\int_{\R^n} |\eta|^{2s} |\hat f|^2
d\eta \right) \left(\int^{\infty}_{0} t^{1-2s} [T^2(t) + T^{'2}(t)] dt\right)  \ ,
\eean
 where $t= |\eta|y$. We next compute the last integral. Multiplying equation (\ref{app.14}) by $T$,
integrating by parts and employing  (\ref{app.18}), we get
\be\la{app.26} \int^{+\infty}_{0} t^{1-2s} [T^2(t) + T^{'2}(t)] dt
=  t^{1-2s}T(t)T'(t)dt \Big|^\infty_0 = \frac{2^{1-2s} \xG(1-s)}{
\xG(s)} \ . \ee

We finally recall the following relation (see, e.g., \cite{FLS}, Lemma 3.1)
\bea\la{app.28}
\int_{\R^n} |\eta|^{2s} |\hat f|^2 d\eta &  = & \frac{c_{n,s}}{2}
 \int_{\R^n} \int_{\R^n}
\frac{|f(x)-f(\xi)|^2}{|x-\xi|^{n+2s}} dx d\xi   \nonumber  \\
& = &
 \frac{s  2^{2s-1} \Gamma(\frac{n+2s}{2})}{ \pi^{\frac{n}{2}}  \Gamma(1-s)}
 \int_{\R^n} \int_{\R^n}
\frac{|f(x)-f(\xi)|^2}{|x-\xi|^{n+2s}} dx d\xi  \ .
\eea
Putting together the last three relations we conclude (\ref{81}).

Finally, taking  into account (\ref{1.20}) we  easily obtain (\ref{82}).

\vspace{2cm}


\begin{thebibliography}{MMMM}

\bibitem [AS]{AS} M. Abramowitz and I. Stegun,  Handbook of Mathematical Functions, with Formulas, Graphs, and
Mathematical Tables.{\em~  NBS Applied Mathematics Series, vol. 55. National Bureau of Standards},
Washington (1964)



\bibitem [ABS]{ABS}
G. Alberti, G. Bouchitt\'e and P. Seppecher, \emph{Phase
Transition with the Line-Tension Effect}, Arch. Rat. Mech. Anal.,
 144 (1), (1998),   1--46.

\bibitem [AFV]{AFV} A. Alvino, A. Ferone and  R. Volpicelli,
 \emph{Sharp Hardy inequalities in the half space with trace remainder term},
 preprint (2011),   arXiv:1105.0335v1.

\bibitem[A]{A} A. Ancona,
\emph{On strong barriers and an inequality of Hardy for domains in
$\Ren$}, J. London Math. Soc. 34 (2), (1986), 274--290.


\bibitem[AK]{AK} D. H. Armitage and U. Kuran,  \emph{The convexity and the superharmonicity of the signed distance
function}, Proc. Amer. Math. Soc. 93 (4), (1985), 598--600.

\bibitem[BFT]{BFT} G. Barbatis, S. Filippas and A. Tertikas,
\emph{A unified approach to improved $L^p$ Hardy inequalities with
best constants}, Trans. Amer. Math. Soc.,  356 (6), (2004),
2169--2196.



\bibitem[BBC]{BBC} K. Bogdan, K. Burdzy and Z-Q. Chen, \emph{Censored
stable processes}, Prob. theory related fields 127 (2003),
89--152.


\bibitem[BD]{BD} K. Bogdan and B. Dyda, \emph{The best constant in a fractional Hardy inequality},
 Math. Nachrichten,   284 (5-6), (2011),  629--638.

\bibitem[BBM]{BBM}
J. Bourgain, H. Brezis and P. Mironescu , \emph{Limiting embedding
theorems
 for $W^{s,p}$  when $s\uparrow 1$ and applications},  J. Anal. Math.  87  (2002), 77--101.


\bibitem[BM]{BM} H. Brezis and M. Marcus,
\emph{Hardy's inequalities revisited}, Dedicated to Ennio De
Giorgi. Ann. Scuola Norm. Sup. Pisa Cl. Sci. {4}, vol. 25 (1997),
217--237.

\bibitem[CC]{CC}
 X. Cabr\'e  and E.  Cinti,
\emph{Energy estimates and 1-D symmetry for nonlinear equations
involving the half-Laplacian}, Discrete Contin. Dyn. Syst. 28 (3),
 (2010), 1179--1206.





\bibitem[CT]{CT} X. Cabr\'e and J. Tan, \emph{ Positive solutions of
nonlinear problems involving the square root of the Laplacian},
Adv. Math.  224 (5),  (2010), 2052--2093.



\bibitem[CS]{CS} L. Caffarelli and L. Silvestre, \emph{ An extension
problem related to the fractional Laplacian}, Comm. P.D.E. 32,
(2007), 1245--1260.

\bibitem[CG]{CG}
Chang S.--Y. A. and Gonzalez M. , \emph{Fractional Laplacian in
conformal geometry}, Adv.  Math.,  226 (2), (2011),   1410--1432.



\bibitem[CKS1]{CKS1}
Z.-Q. Chen, P. Kim and R. Song, \emph{Two-sided heat kernel
estimates for censored stable-like processes},
  Probab. Theory Related Fields  146 (3-4),  (2010),   361--399.

\bibitem[CKS2]{CKS2}
Z.-Q. Chen, P. Kim and R. Song,  \emph{Heat Kernel Estimates for
Dirichlet Fractional Laplacian}, J. Eur. Math. Soc. 12(5), (2010),
1307--1329.


\bibitem[CSo]{CSo} Z-Q. Chen and R. Song, \emph{Hardy inequality for censored stable
processess},
 Tohoku Math. J.,  55 (2), (2003), 439--450.

\bibitem[D1]{D1} E. B. Davies,
\emph{The Hardy constant}, Quart. J. Math.,  46 (2), (1995),
417--431.


\bibitem[D2]{D2}  E. B. Davies, \emph{A review of Hardy inequalities}, Oper. Theory Adv.
Appl.,  110, (1999), 55--67.

\bibitem[DDM]{DDM} J. Davila, L. Dupaigne and M. Montenegro, \emph{The
extremal solution of a boundary reaction problem},
 Commun. Pure Appl. Anal.  7 (4),  (2008),   795--817.

\bibitem[D]{D}  B. Dyda,
 \emph{Fractional Hardy-Sobolev-Maz'ya inequality on balls and halfspaces}, preprint (2010)
 arXiv:1004.5146.

\bibitem [E]{E} L. C. Evans,  \emph{Partial Differential Equations},  American Mathematical Society
Providence, RI (1998)


 \bibitem[FMT]{FMT2} S. Filippas, V. Maz'ya and A. Tertikas,
\emph{Critical Hardy-Sobolev Inequalities}, J. Math. Pures Appl.
(9), 87 (1),
 (2007), 37-56.

\bibitem[FS]{FS} R. L. Frank and R. Seiringer, \emph{Sharp fractional Hardy inequalities in
half-spaces},
 In: Around the research of Vladimir Maz'ya, A. Laptev (ed.), 161 - 167,
 International Mathematical Series 11 (2010).

\bibitem[FLS]{FLS} R. L. Frank, E. H. Lieb and R. Seiringer, \emph{Hardy-Lieb-Thirring inequalities
 for fractional Schroedinger operators}, J. Amer. Math. Soc. 21 (4), (2008), 925--950.


\bibitem[FL]{FL}  R. L. Frank and M. Loss,
\emph{Hardy-Sobolev-Maz'ya inequalities for arbitrary domains},
preprint (2011), arXiv:1102.4394v1.

 \bibitem[Gk]{Gk} K. Gkikas, \emph{Hardy and Hardy-Sobolev inequalities and their Applications}, Ph.D. thesis 2011.

\bibitem[G]{G} M. Gonzalez, \emph{Gamma convergence of an energy functional related to the
 fractional Laplacian},   Calc. Var. Partial Differential Equations,  36 (2), (2009),  173--210.



\bibitem[KK]{KK} J. Kinnunen  and R. Korte, \emph{Characterizations for the Hardy Inequality},
  Around the research of Vladimir Maz'ya, A. Laptev (ed.),
239--254,
 International Mathematical Series 11 (2010).

\bibitem[L]{L} Leoni G. \emph{A first course in Sobolev Spaces}, Grad. Studies in Math.,
American Math. Society, (2009).

\bibitem[LLL]{LLL}  R. T. Lewis, J. Li and Y. Y. Li,
 \emph{A geometric characterization of a sharp Hardy inequality}, preprint (2011), arXiv:1103.5429.


\bibitem[LS]{LS}
M. Loss and C.  Sloane, \emph{Hardy inequalities for fractional
integrals on general domains}, J. Funct. Anal. 259 (6), (2010),
1369--1379.



\bibitem[M]{Maz}  V. Maz'ya,
\emph{Sobolev spaces with applications to elliptic partial
differential equations},
 Second, revised and augmented edition.
  Springer, Heidelberg, 2011.

\bibitem[MS]{MS}
V. Maz'ya and T. Shaposhnikova,  \emph{On the Bourgain, Brezis,
and Mironescu theorem concerning limiting
 embeddings of fractional Sobolev spaces},  J. Funct. Anal.  195 (2)  (2002),  230--238; Erratum  J. Funct. Anal.,
201 (1),  (2003),  298--300.

\bibitem[N]{N} H. M. Nguyen,
 \emph{$\Gamma$--convergence, Sobolev norms and BV functions},
Duke Math. J., 157 (3), (2011), 495--533.

\bibitem[PSV]{PSV}
G. Palatucci, O. Savin and E. Valdinoci, \emph{Local and global
minimizers for a
 variational energy involving a fractional norm}, preprint (2011),   arXiv:1104.1725.

\bibitem[P]{P} G. Psaradakis, \emph{$L^1$ weighted Hardy inequalities}, preprint  (2011).

\bibitem[SV1]{SV1}
O. Savin and E. Valdinoci, \emph{Density estimates for a
variational model driven by the Gagliardo norm}, preprint (2010),
arXiv:1007.2114.


\bibitem[SV2]{SV2}
O. Savin and E. Valdinoci, \emph{$\Gamma$-convergence for nonlocal
phase transitions},  preprint (2010),   arXiv:1007.1725.


\bibitem[S]{S}
  C. A. Sloane,  \emph{A Fractional Hardy-Sobolev-Maz'ya Inequality on the Upper
  Halfspace},
preprint (2010), arXiv:1004.4828.

\bibitem[St]{St} E. M. Stein, \emph{Singular Integrals and Differentiability Properties of Functions},
Princeton Univ. Press, 1970.

\bibitem[T]{T}  J. Tan,
\emph{The Brezis  Nirenberg type problem involving the square root
of the Laplacian},
 Calc. Var. Partial Differential Equations, 42 (1-2), (2011),
 21--41.



\end{thebibliography}
\end{document}